\documentclass[11pt,letterpaper]{amsart}
\usepackage{graphicx} 
\usepackage{amsthm}
\usepackage{amsfonts}
\usepackage{amssymb}
\usepackage{amsmath}
\usepackage{mathtools}
\usepackage[maxbibnames=99]{biblatex}
\usepackage{xcolor}
\usepackage{esint} 
\usepackage[citecolor=blue,colorlinks]{hyperref}

\newtheorem{innercustomgeneric}{\customgenericname}
\providecommand{\customgenericname}{}
\newcommand{\newcustomtheorem}[2]{%
  \newenvironment{#1}[1]
  {%
   \renewcommand\customgenericname{#2}%
   \renewcommand\theinnercustomgeneric{##1}%
   \innercustomgeneric
  }
  {\endinnercustomgeneric}
}

\newcustomtheorem{customthm}{Theorem}
\newcustomtheorem{customlemma}{Lemma}
\newcustomtheorem{customcorollary}{Corollary}
\newcustomtheorem{customremark}{Remark}

\newtheorem{thm}{Theorem}[section]
\newtheorem{corollary}[thm]{Corollary}
\newtheorem{lemma}[thm]{Lemma}

\newtheorem*{thm*}{Theorem}
\newtheorem*{corollary*}{Corollary}
\newtheorem*{lemma*}{Lemma}
\newtheorem*{proposition*}{Proposition}

\theoremstyle{definition}
\newtheorem{definition}[thm]{Definition}
\newtheorem*{definition*}{Definition}
\newtheorem{remark}[thm]{Remark}

\newtheorem*{remark*}{Remark}

\newcommand{\bb}[1]{\mathbb{#1}}
\newcommand{\ssf}[1]{\mathsf{#1}}

\newcommand{\Lip}{\mathsf{Lip}}

\newcommand{\Ric}{\mathsf{Ric}}

\newcommand{\Hess}{\mathsf{Hess}}

\newcommand{\aH}{\mathsf{H}}

\newcommand{\X}{\mathsf{X}}
\newcommand{\Vol}{\mathsf{Vol}}

\newcommand{\sd}{\mathsf{d}}
\newcommand{\sL}{\mathsf{L}}
\newcommand{\de}{\mathrm{d}}

\makeatletter
\newsavebox\myboxA
\newsavebox\myboxB
\newlength\mylenA

\newcommand*\xoverline[2][0.75]{%
    \sbox{\myboxA}{$\m@th#2$}%
    \setbox\myboxB\null
    \ht\myboxB=\ht\myboxA%
    \dp\myboxB=\dp\myboxA%
    \wd\myboxB=#1\wd\myboxA
    \sbox\myboxB{$\m@th\overline{\copy\myboxB}$}
    \setlength\mylenA{\the\wd\myboxA}
    \addtolength\mylenA{-\the\wd\myboxB}%
    \ifdim\wd\myboxB<\wd\myboxA%
       \rlap{\hskip 0.5\mylenA\usebox\myboxB}{\usebox\myboxA}%
    \else
        \hskip -0.5\mylenA\rlap{\usebox\myboxA}{\hskip 0.5\mylenA\usebox\myboxB}%
    \fi}
\makeatother

\title[Splitting  for manifolds with a convex boundary component]{A splitting theorem for manifolds with a convex boundary component and applications}
\author{Alessandro Cucinotta  \and Andrea Mondino}

\begin{document}

\begin{abstract}
   We prove a warped product splitting theorem for manifolds with Ricci curvature bounded from below in the spirit of [Croke-Kleiner, \emph{Duke Math.\;J}.\;(1992)], but instead of asking that one boundary component is compact and mean-convex, we require that it is parabolic and convex. 
    We then deduce several applications, including splitting theorems and first Betti number rigidity results for
    \begin{itemize}
    \item  $3$-manifolds with non-negative Ricci curvature, 
    \item $4$-manifolds with weakly bounded geometry, non-negative $2$-Ricci curvature, scalar curvature $\geq 1$.
    \end{itemize}
    In particular, the latter aswers to a rigidity question posed by [Chodosh-Li-Stryker, \emph{JEMS},\;(2024)].
   The proofs rely on a metric gluing of Riemannian manifolds with boundary, resulting in a non-smooth metric space.
To address this lack of smoothness, we employ synthetic tools specifically developed for non-smooth settings, with a focus on those based on optimal transportation.
\end{abstract}

\maketitle

\section{Introduction}

The celebrated Cheeger-Gromoll splitting theorem \cite{SplittingCheegerGromoll} states that a (complete, connected) Riemannian manifold with non-negative Ricci curvature containing a line is isometrically a product, where one of the factors is a real line. 
A variant of this result was proved by Kasue in \cite{Kasue}. There, it is shown that
a Riemannian manifold $(M,g)$ with mean-convex compact boundary $\partial M$ and non-negative Ricci curvature, containing a half line
with initial point in $\partial M$, again is isometrically a product. In the latter case, the splitting factor is a half line.
This result fails without the compactness assumption on the boundary, as one sees considering the epigraph of a strictly convex function in Euclidean space.
In the same work, Kasue shows that a manifold with multiple mean-convex boundary components, one of which is compact, and non-negative Ricci curvature is also a product, where one of the factors is a bounded interval of the real line. 
A generalization of Kasue's results to manifolds with  Ricci curvature bounded below by a \emph{negative} constant was given by Croke and Kleiner in \cite{warpedsplit}, and is recalled below.

Let us first fix some notation. In a Riemannian manifold $(M^n,g)$ with boundary $\partial M$ and external unit normal $\nu$, we denote by $\mathrm{I\!I}^M_{\partial M}$ the second fundamental form of the boundary: $$\mathrm{I\!I}^M_{\partial M}(X,Y):=g(\nabla_X \nu,Y), \quad \text{for all } X,Y\in T_xM,\; x\in \partial M.$$
We set 
$$H^M_{\partial M}:=\frac{1}{n-1}\, \mathrm{tr}(\mathrm{I\!I}^M_{\partial M})$$ 
to be the associated mean curvature (the upper script $M$ will be omitted if there is no ambiguity).
We say that a manifold with boundary $(M,g)$ is convex (resp.\ mean convex) if its boundary satisfies $\mathrm{I\!I}_{\partial M} \geq 0$ (resp.\ $H_{\partial M} \geq 0$). Similarly, a connected component $S$ of $\partial M$ is said to be convex (resp.\ mean convex) if it satisfies $\mathrm{I\!I}_{S} \geq 0$ (resp. $H_{S} \geq 0$). The intrinsic Ricci curvature of $S$ is denoted by $\Ric_{S}$.

\begin{thm*}[{\cite[Theorem 1]{warpedsplit}}]
    Let $\delta \in \{0,1\}$ and let $(M^n,g)$ be a Riemannian manifold with $\Ric_{M} \geq -(n-1) \delta$ and whose boundary consists of two connected components $\partial M=S_1 \cup S_2$, one of which is compact.  Assume that the mean curvature of $S_1$ is $\geq \delta$ and the mean curvature of $S_2$ is $\geq -\delta$. 
    Then, there exists $l>0$ such that $M$ is isometric to $S_1 \times [0,l]$ with
    the warped product metric $\de s^2 = e^{-2\delta t}g_1 + \de t^2$, where $g_1$ is the metric on $S_1$. Moreover, $\Ric_{S_1} \geq 0$.
\end{thm*}

In Theorems \ref{CT1}, \ref{CT1*}, and \ref{CT2} below, we prove variants of this result. On the one hand, instead of asking for the compactness of one of the boundary components, we ask for its parabolicity, which is a weaker condition. For example, $\bb{R}^2$ is parabolic. On the other hand, we assume that the parabolic boundary component has a lower bound on the second fundamental form, rather than on the mean curvature, and that it has itself nonnegative Ricci curvature. Whether Theorems \ref{CT1}, \ref{CT1*}  and \ref{CT2} hold with just a mean curvature lower bound on the parabolic boundary component is an open question, which we discuss in Remark \ref{R2} and in Section \ref{S5}.

Recall that a manifold is said to be \emph{parabolic} if it admits no positive fundamental solution for the Laplacian. For example, while $\bb{R}^2$ is parabolic, $\bb{R}^n$ with $n \geq 3$ is not. More generally, a manifold $(M,g)$ with non-negative Ricci curvature is parabolic if and only if $$\int_1^\infty \frac{t}{\ssf{Vol}(B_t(x))} \, \de t= +\infty,$$ for some $x \in M$. 
Next, we state the main results of this paper.

\begin{thm} \label{CT1}
    Let $(M^n,g)$ be a Riemannian manifold with $\Ric_{M} \geq 0$ and disconnected mean-convex boundary. Assume that one connected component $S \subset \partial M$ is parabolic, convex, and satisfies $\Ric_{S} \geq 0$.
    
    Then, there exists $l>0$ such that $M$ is isometric to $S \times [0,l]$ with
    the product metric.
\end{thm}

The next two theorems provide variants of Theorem \ref{CT1} where the Ricci curvature is bounded from below by a \emph{negative} constant (compare with Croke-Kleiner's result above). 

Theorem \ref{CT1*} deals with the case where a parabolic boundary component satisfies $\mathrm{I\!I} \geq 1$, while the other boundary components satisfy $H \geq -1$.

\begin{thm} \label{CT1*}
    Let $(M^n,g)$ be a Riemannian manifold with $\Ric_{M} \geq -(n-1)$ and disconnected boundary with $H_{\partial M} \geq -1$. 
    Assume that one connected component $S \subset \partial M$ is parabolic, has $\mathrm{I\!I}_{S} \geq 1$, and $\Ric_{S} \geq 0$.
    
    Then, there exists $l>0$ such that $M$ is isometric to $S \times [0,l]$ with
    the warped product metric $\de s^2 = e^{-2 t}g_S + \de t^2$, where $g_S$ is the metric on $S$.
\end{thm}

Theorem \ref{CT2} treats the case where a parabolic boundary component satisfies $\mathrm{I\!I} \geq -1$, while the other boundary components satisfy $H \geq 1$.

\begin{thm} \label{CT2}
     Let $(M^n,g)$ be a Riemannian manifold, with $\Ric_{M} \geq -(n-1)$ and disconnected boundary. Assume that one connected component $S \subset \partial M$ is parabolic, has $\mathrm{I\!I}_S \geq -1$, and $\Ric_S \geq 0$. Assume also that $H_{\partial M \setminus S} \geq 1$.
 
    Then, $S_1:=\partial M \setminus S$ is connected, and there exists $l>0$ such that $M$ is isometric to $S_1 \times [0,l]$ with
    the warped product metric $\de s^2 = e^{-2t}g_1 + \de t^2$, where $g_1$ is the metric on $S_1$.
    \end{thm}



\begin{remark}[On the parabolicity assumption]
In Theorem \ref{CT1}, the parabolicity assumption cannot be omitted. To see this, consider the manifold $M \subset \mathbb{R}^4$ given by the region bounded by a catenoid (see \cite{Catenoid1, Catenoid2}) and a disjoint hyperplane. In this setting, $M$ satisfies all the hypotheses of Theorem \ref{CT1}, except for the parabolicity of the convex boundary component, since $\mathbb{R}^3$ is not parabolic. However, the conclusion of the theorem does not hold. Similar counterexamples are expected to exist for Theorems \ref{CT1*} and \ref{CT2}. This phenomenon is related to the study of the Half Space Property (see \cite{CMMR} and the discussion below in the introduction) in the context of warped products.
\end{remark}

\begin{remark}[On the curvature assumptions for the parabolic boundary component] \label{R2}
Given the result by Croke-Kleiner, it is natural to ask whether it is possible to replace the conditions
\begin{itemize}
\item $\mathrm{I\!I}_{S} \geq 0$, $\Ric_S \geq 0$ in Theorem \ref{CT1};
\item $\mathrm{I\!I}_{S} \geq 1$, $\Ric_S \geq 0$ in Theorem \ref{CT1*};
\item $\mathrm{I\!I}_{S} \geq -1$, $\Ric_{S} \geq 0$ in Theorem \ref{CT2};
\end{itemize}
with the considerably weaker conditions
\begin{itemize}
\item $H_{S} \geq 0$ in Theorem \ref{CT1};
\item $H_{S} \geq 1$ in Theorem \ref{CT1*};
\item $H_{S} \geq -1$ in Theorem \ref{CT2}.
\end{itemize}
It is likely that the answer to this question is affirmative (see Section \ref{S5} for a partial answer). 

The stronger boundary assumptions in the aforementioned theorems stem from the gluing construction underlying our proof, which produces a sufficiently regular metric space only under appropriate control of the second fundamental form of the boundary. The following example provides a useful illustration of the role played by boundary assumptions in gluing constructions. By \cite{schlichting2012gluing}, if $M$ has $\Ric \geq 0$ and $\mathrm{I\!I}_{\partial M} \geq 0$, then its double (i.e. the gluing of a copy of $M$ to itself along the boundary) is a $C^0$ manifold whose metric can be approximated locally by metrics with almost nonnegative Ricci curvature. On the other hand, if $M$ has $\Ric \geq 0$ and $H_{\partial M} \geq 0$, then the metric of its double can only be approximated locally by metrics with almost nonnegative \emph{scalar} curvature, which is a much weaker condition. For instance, using the surgery results in \cite{GromovLawson1,GromovLawson2}, one can construct manifolds with nonnegative scalar curvature and two totally geodesic boundary components which are not products.
\end{remark}

\begin{remark}[How to construct manifolds as in Theorems \ref{CT1}, \ref{CT1*}, \ref{CT2}]
Using the extensive literature on the existence and regularity of minimal surfaces, one
can construct manifolds satisfying the assumptions of Theorem \ref{CT1} starting from fairly general assumptions.  Indeed, let $(N,g)$ be a smooth Riemannian manifold where any \emph{stable} minimal hypersurface $\Sigma \subset N$ is totally geodesic, parabolic, and has $\Ric_\Sigma \geq 0$ (existence of such manifolds follows, for instance, from  \cite{SchoenYauStability} and \cite{chodosh2024completestableminimalhypersurfaces}). Then, any smooth region $M \subset N$ bounded by the disjoint union of a minimal hypersurface and a \emph{stable} minimal hypersurface satisfies the boundary assumptions of Theorem \ref{CT1}. Such regions can be constructed either assuming that $N$ has disconnected minimal boundary using a variation of \cite[Lemma $4$, Section $11$]{MeeksSimonYau} (this is the strategy used in Corollaries \ref{CC2} and \ref{CC3:intro} below), or assuming that $N$ supports a free closed action of $\bb{Z}$ using an argument from \cite{SchoenYauStability} (this is the strategy used in Theorem \ref{T|Anderson} below). Similar constructions are likely to be possible for manifolds as in Theorems \ref{CT1*} and \ref{CT2}, replacing the use of minimal surfaces with the use of $\mu$-bubbles (see \cite{zGromovbubbles1,GromovLecturesScalar}).
\end{remark}

We next discuss the main applications. As an application of Theorem \ref{CT1}, we deduce splitting results for low dimensional manifolds with suitable curvature assumptions and disconnected mean-convex boundary, with no further assumptions (not even parabolicity) on the boundary. We stress that these results, in particular, do not follow from Kasue's Theorem.

The first such result is Corollary \ref{CC2}, which shows that a $3$-manifold with non-negative Ricci curvature and disconnected mean-convex boundary is a product. To the best of our knowledge, this statement is original. The proof relies on a new criterion to determine whether a manifold is a product (see Theorem \ref{CT4}), on the results from \cite{SchoenYauStability}, and on \cite[Lemma $4$, Section $11$]{MeeksSimonYau}.

\begin{corollary}\label{CC2}
    Let $(M^3,g)$ be a Riemannian manifold with non-negative Ricci curvature and disconnected mean-convex boundary. Then, $(M,g)$ splits isometrically as $\Sigma \times [0,l]$, for some manifold $(\Sigma,g')$ and some $l>0$.
\end{corollary}

\begin{remark}[On a related result by Anderson-Rodriguez] \label{R1}
In \cite[Theorem 3]{AndersonRodriguez}, Anderson and Rodriguez proved a result similar to Corollary \ref{CC2}, assuming an additional \emph{uniform upper bound on the sectional curvature} and that the disconnected boundary is \emph{minimal}. 
\end{remark}

As a second application of Theorem \ref{CT1}, combined with  \cite[Theorem 1.10]{chodosh2024completestableminimalhypersurfaces} and the tools used to prove Corollary \ref{CC2} (i.e.,\ Theorem \ref{CT4}), we deduce a splitting result for $4$-manifolds with non-negative $2$-Ricci curvature, scalar curvature greater than $1$,  weakly bounded geometry  and mean-convex disconnected boundary (cf.\;\cite[Theorem 5.2]{FrankelMaximum} for the case of  non-negative \emph{sectional} curvature and  \emph{minimal} boundary). 
We refer to Definitions \ref{D2} and \ref{D1} (and the subsequent comments) for the precise definitions of manifolds with weakly bounded geometry and  non-negative $2$-Ricci curvature (denoted $\Ric_2 \geq 0$), and for their geometric relevance. We only mention here that $\Ric_{2} \geq 0$ is an intermediate condition between $\mathsf{Sec} \geq 0$ and $\Ric \geq 0$, in the sense that $\mathsf{Sec} \geq 0 \Rightarrow \Ric_{2} \geq 0 \Rightarrow \Ric \geq 0$.

\begin{corollary}\label{CC3:intro}
    Let $(M^4,g)$ be a Riemannian manifold with $\Ric_2 \geq 0$, scalar curvature $\geq 1$, and weakly bounded geometry. Let $N^4 \subset M^4$ be a smooth submanifold with mean-convex disconnected boundary. Then, $N=\Sigma \times [0,l]$ isometrically, for some manifold $(\Sigma,g')$ and $l>0$.
\end{corollary}

Through a classical argument linking the non-vanishing of the first Betti number to existence of minimal hypersurfaces, Corollaries \ref{CC2} and \ref{CC3:intro} can be used to obtain the following result.

\begin{thm} \label{T|Anderson}
    Let $(M,g)$ be a manifold without boundary such that either
    \begin{itemize}
    \item  $n=3$, and $\Ric \geq 0$,
    \end{itemize}
    or
    \begin{itemize}
    \item $n=4$, $\Ric_2 \geq 0$, $\mathrm{Scal} \geq 1$, and $M$ has weakly bounded geometry.
    \end{itemize}
    If $b_1(M) \neq 0$, then the universal covering of $M$ splits a line isometrically.
\end{thm}
In the case of $3$-manifolds with nonnegative Ricci curvature, the previous theorem was known from \cite{GangLiu} (see also \cite{cucsemmagn}), with a different proof. On the other hand, the result is new for $4$-manifolds with $\Ric_2 \geq 0$, scalar curvature $\geq 1$, and weakly bounded geometry, and it answers a rigidity question from \cite[Comments after Theorem 1.13]{chodosh2024completestableminimalhypersurfaces}.

\medskip

As an application of Theorems \ref{CT1*} and \ref{CT2}, we obtain a slice theorem for warped products over parabolic manifolds with non-negative Ricci curvature, see Corollary \ref{CC1} below.
Let us give some context first.
Following \cite{RSS}, a manifold $(M,g)$ is said to have the Half Space Property if the only (properly embedded) minimal hypersurfaces of $M \times \bb{R}$ contained in a half-space are the horizontal slices $M \times \{t\}$. In recent years, several results have been obtained by different authors (see, among others, \cite{DBH,HLRS,EH,DBM,CMMR, DingCap}). 
Corollary \ref{CC1} provides a result in the spirit of the Half Space Property for (possibly warped) products. The main differences with the aforementioned classical half-space results are that:
\begin{itemize}
\item Corollary \ref{CC1} holds for mean-convex boundaries in product manifolds, not only for minimal hypersurfaces;
\item we obtain also a half-space result for sets whose boundary has mean curvature bounded below by $1$ (or $-1$, depending whether they lie in the lower or upper half space), in warped products with negative curvature.
\end{itemize}

\begin{corollary} \label{CC1}
    Let $\delta \in \{0,1\}$. Let $(M^n,g)$ be a parabolic manifold with $\Ric_M \geq 0$ and let $M \times \bb{R}$ be equipped with the metric $\de s^2 = e^{-2 t \delta}g + \de t^2$. If $E \subset M \times (0,+\infty)$ is a smooth closed set with connected boundary and mean curvature $H_{\partial E} \leq \delta$, then $E=M \times [a,+\infty)$, for some $a > 0$.
    If $E \subset M \times (-\infty ,0)$ is a smooth closed set with connected boundary and mean curvature $H_{\partial E} \leq -\delta$, then $E=M \times (-\infty,a]$, for some $a < 0$.
\end{corollary}


\begin{remark}[On a related result by Montiel]
    In \cite{Montiel}, Montiel proved that if $(M,g)$ is compact with non-negative Ricci curvature, and $M \times \bb{R}$ is equipped with the metric $\de s^2=e^{-2t}g+\de t^2$, then any hypersurface $\Sigma \subset M \times \bb{R}$ of constant mean curvature, that is locally a graph on $M$, must be a slice. 
    Some generalizations of this result later appeared in \cite{HalfSpaceWarped,HalfSpaceSlab,Cavalcante}. Although similar in spirit, Montiel's theorem is independent of Corollary \ref{CC1}: indeed, it requires the mean curvature to be constant while we require the inequality $H\leq \pm 1$ and we do not assume the hypersurface to be locally a graph.    
\end{remark}

\medskip

We briefly outline the strategy of the proof of Theorem \ref{CT1}.
Consider the metric space $(\X,\sd)$ obtained by gluing $M$ and $S \times [0,+\infty)$ along their isometric boundaries.  Let us emphasize that $(\X, \sd)$ is a non-smooth space: the glued Riemannian metric is only of class $C^{0,1}$ along $S$. This lack of smoothness motivates the use of synthetic tools developed for non-smooth spaces, particularly those based on optimal transport.

Consider the distance function $\sd_{S_2}$ from $S_2:=\partial M \setminus S$ in the glued space and its restriction to $S \times [0,+\infty)$. We denote by $\sd_{S}$ the distance from $S$ in $S \times [0,+\infty)$. If we can show that $\sd_{S_2}$ is constant on $S \times \{0\}$, then the statement follows by standard arguments.

Thanks to the assumption on the fundamental form of $S$ (and previous results on the Laplacian of distance functions in metric measure spaces with synthetic Ricci curvature lower bounds, see  \cite{Giglimem, CM18,  ketterer2023rigidity, Weak, lapb}), it holds $\Delta \sd_{S_2} \leq 0$ on $S \times (0,+\infty)$ in distributional sense. 
\\In particular, it holds $\Delta (\sd_{S_2}-\sd_{S}) \leq 0$ and, calling $\nu$ the exterior normal of $S \times \{0\}$ in $S \times [0,+\infty)$, it also (formally) holds $\nabla (\sd_{S_2}-\sd_{S}) \cdot \nu \geq 0$ on $S \times \{0\}$. 
If $S \times [0,+\infty)$ has sufficiently small volume growth at infinity, we then deduce that $\sd_{S_2}-\sd_{S}$ is constant by an integration by parts argument, cf.\;\cite{ParabolicitaPigola}. 
Since in general $S \times [0,+\infty)$ does not satisfy the required volume growth at infinity, we multiply the volume measure of $S \times [0,+\infty)$ by a suitable weight, cf.\;\cite{C1}. 
Making an appropriate choice of such a weight, we obtain that $\sd_{S_2} - \sd_{S}$  satisfies the previous Laplacian bound also in the weighted space, which, in addition, has the desired volume growth. We then carry out a weighted integration by parts argument to obtain that $\sd_{S_2} - \sd_{S}$ is constant. Hence, $\sd_{S_2}$ is constant on $S \times \{0\}$.

\medskip
\noindent
\textbf{Acknowledgements}. The authors wish to thank Mattia Magnabosco for useful discussions. They also wish to thank the reviewers for their careful reading of the manuscript and their useful suggestions. 

A.\,M.\;acknowledges support from the European Research Council (ERC) under the European Union’s Horizon 2020 research and innovation programme, grant agreement No. 802689 “CURVATURE”.

\section{Preliminary lemmas}

We first introduce a few optimal transport tools that will be used in the proof of Lemma \ref{l|Laplacian in the glued space}. 

If Lemmas \ref{l|Laplacian in the glued space} and \ref{CL2} -- whose statements do not rely on optimal transport theory -- are taken as given, this section may be skipped, and one can proceed directly to the proofs of Theorems \ref{CT1}, \ref{CT1*}, and \ref{CT2}.

\medskip
Let $\X$ be an $n$-dimensional smooth manifold (possibly with boundary) equipped with a locally Lipschitz continuous metric $g$. The metric $g$ induces a distance $\sd$ and a volume measure, which coincides with the $n$-dimensional Hausdorff measure $\aH^n$ induced by $\sd$. 
Given an open set $\Omega \subset \X$, we denote by $\Lip(\Omega)$ and $\Lip_c(\Omega)$ respectively Lipschitz continuous and compactly supported Lipschitz continuous functions on $\Omega$. 

Assume that the metric space $(\X,\sd)$ is complete. On the metric measure space $(\X,\sd,\aH^n)$ one can define the Wasserstein distance (w.r.t.\;to the distance squared cost) between two probability measures with finite second moment. The set of probability measures with finite second moment is denoted $\mathcal{P}_2(\X)$. 
Given $N \in (1,+\infty)$ and $\mu=\rho \aH^n \in \mathcal{P}_2(\X)$,  the R\'enyi entropy of $\mu$ with respect to $\aH^n$ is defined as
\[
U_{N}(\mu|\aH^n):=-\int \rho^{1-1/N} \,\de \aH^n.
\]
Given $K \in \bb{R}$, two measures $\mu_1=\rho_1\aH^n,\mu_2=\rho_2 \aH^n \in \mathcal{P}_2(\X)$, and an optimal plan $\Pi$ between them, we define
\begin{align*}
T^{(t)}_{K,N}(\Pi|\aH^n):=\int_{\X \times \X} [ &\tau^{1-t}_{K,N}(\sd(x_1,x_2))\rho_1^{-1/N}(x_1) \\ & \quad +
\tau^{t}_{K,N}(\sd(x_1,x_2))\rho_2^{-1/N}(x_2)] \, d\Pi (x_1,x_2),
\end{align*}
where, for every $s \in [0,1]$, $\tau^{s}_{\cdot,\cdot}(\cdot)$ is an appropriate distortion coefficient (see \cite{S2}).

We say that \emph{$\X$ satisfies the ${\rm CD}(K,N)$ condition in $U \subset \X$}, if the following holds.
For every pair $\mu_1=\rho_1\aH^n,\mu_2=\rho_2 \aH^n \in \mathcal{P}_2(\X)$ of measures supported in $U$, there exists a Wasserstein geodesic $\{\xi_t\}_{t \in [0,1]} \subset \mathcal{P}_2(\X)$ from $\mu_1$ to $\mu_2$ and an optimal coupling $\Pi$ of $\mu_1$ and $\mu_2$ such that, for every $t \in [0,1]$ and every $N' \geq N$, it holds
\[
U_{N'}(\xi_t|\aH^n) \leq -T^{(t)}_{K,N'}(\Pi|\aH^n).
\]
This condition differs from the $\mathrm{CD}_{loc}(K,N)$ condition that appears in the literature (see for instance \cite{BacherSturm, CavallettiMilman}), since we are not requiring that the convexity of the entropy holds in a neighbourhood of every point. For more background on curvature dimension conditions in metric (measure) spaces, we refer to the foundational works \cite{S1,S2,V}.

\medskip
We now turn our attention to Laplacians of functions on $\X$, referring to \cite{Giglimem, CM18,  ketterer2023rigidity, Weak, lapb} for more details and results on Laplacians in metric measure spaces with synthetic Ricci lower bounds. Let $\Omega \subset \X$ be an open set.
We say that a Radon measure $\mu$ on $\Omega$ is the distributional Laplacian of $f \in \Lip(\Omega)$ (and we write $\Delta f=\mu$) if, for every $\phi \in \Lip_c(\Omega)$, it holds
\[
-\int_\Omega \nabla f \cdot \nabla \phi \,\de \aH^n= \int \phi \, d\mu.
\]
We remark that since we are working on a smooth manifold with a continuous metric, the product $\nabla f \cdot \nabla \phi$ is well defined $\aH^n$-almost everywhere. In particular, for our purposes, we do not need to consider more general notions of gradients on metric spaces.

\medskip
We now recall some facts about the localization technique for $\textrm{CD}(K,N)$ spaces. We refer to \cite{MondinoInventiones, CM18} for the proofs, as well as for the definitions of transport set of a Lipschitz function, transport rays and disintegration of a measure.

Let $N > 1$ and $K \in \bb{R}$. Assume that $\X$ is $\textrm{CD}(K,N)$ in a neighbourhood of every point and that it is non-branching (i.e.\;geodesics do not branch). Let $\phi \in \Lip(\X)$ be a $1$-Lipschitz function and consider the partition of its transport set into transport rays $\{r_\alpha\}_{\alpha \in Q}$, $Q$ being a set of indices with a measure $\ssf{q}$ induced by the partition.
Consider the associated disintegration of the measure $\aH^n$ (restricted to the transport set of $\phi$) into measures $\{h_\alpha\}_{\alpha \in Q}$, each concentrated on the corresponding transport ray. 

Then, $\ssf{q}$-almost every measure $h_\alpha$ is absolutely continuous w.r.t.\;the Lebesgue measure of the corresponding transport ray, it admits a locally Lipschitz continuous density (that we still denote by $h_\alpha$), and it satisfies in weak sense 
\[
    (h_\alpha^{\frac{1}{N-1}})''+\frac{K}{N-1} h_\alpha^{\frac{1}{N-1}} \leq 0.
\]
It was shown in \cite{CM18} that, if the function $\phi$ is the distance function from a closed set $E \subset \X$, we also have
\[
\Delta \phi = (\log h_\alpha)'+[\Delta \phi]^{sing} \quad \text{in } \X \setminus E \text{ in distributional sense}.
\]
In the previous equation, $[\Delta \phi]^{sing}$ is a negative singular measure (w.r.t. $\aH^n$), while the derivative $(\log h_\alpha)'$ is taken orienting the transport rays from $E$ to $\X \setminus E$.

The next Lemma deals with the Laplacian of distance functions in metric spaces resulting from the gluing of two manifolds along their isometric boundaries. In smooth manifolds, a lower bound on the Ricci curvature gives upper bounds on the Laplacian of distance functions. Analogously, in the lemma below, we need the resulting metric space to have Ricci curvature bounded from below in a synthetic sense (i.e.\ to be a $\mathrm{CD}(K,N)$ space) on the region where the distance function is strictly positive. To achieve this, we need a lower Ricci curvature bound on the glued manifolds, and that the sum of the second fundamental forms of the glued boundaries is nonnegative (compare with Remark \ref{R2} from the Introduction).
\begin{lemma}
    \label{l|Laplacian in the glued space}
    Let $\delta \in \{0,1\}$ and let $(M^n,g)$ be a Riemannian manifold with $\Ric_{M} \geq -(n-1) \delta$. Let $S_1$ be a union of connected components of $\partial M$. Assume that $\partial M \setminus S_1 \neq \emptyset$ and set $S_2:=\partial M \setminus S_1$. 
    Let $(N,\tilde{g})$ be a second manifold (possibly disconnected) with $\Ric_{N} \geq -(n-1) \delta$, and whose boundary $\partial N$ is isometric to $S_1$.
    Let $(\X,\sd)$ be the metric space obtained by gluing $M$ and $N$ along the two isometric copies of $S_1$, and let $\sd_{S_2}$ be the distance from $S_2$ in $\X$. Assume that $\mathrm{I\!I}^M_{S_1}+\mathrm{I\!I}^N_{\partial N} \geq 0$.
    \begin{enumerate}
    \item If $H_{S_2} \geq -\delta$, then $\Delta \sd_{S_2} \leq (n-1) \delta$ in $\X \setminus S_2$ in distributional sense.
    \item If $H_{S_2} \geq \delta$, then $\Delta \sd_{S_2} \leq -(n-1) \delta$ in $\X \setminus S_2$ in distributional sense.
    \item Let $S$ be a connected component of $S_1$, and let $N=N_S \cup N_{S_1 \setminus S}$, where $N_S$ and $N_{S_1 \setminus S}$ are disjoint manifolds, and $\partial N_S \cong S$, $\partial N_{S_1 \setminus S} \cong S_1 \setminus S$. Let $D \geq 0$ be the distance between $S$ and $S_2$ in $\X$. If $H_{S} \geq -\delta$ in $M$, then $\Delta \sd_{S} \leq (n-1) \delta$ in the $D$-neighbourhood of $S$ in $\X \setminus N_S$ in distributional sense.
    \end{enumerate}
        \begin{proof}
        From now on, we refer to $S_1$, $S_2$, $M$ and $N$ as subsets of $\X$. Balls $B_r(x) \subset \X$ are always defined w.r.t. the distance $\sd$ on $\X$. The space $\X$ admits a smooth structure such that
        the metric $g^\X$ on $\X$ defined by $g^\X_{|M}=g$ and $g^\X_{|N}=\tilde{g}$ is locally Lipschitz continuous and it induces the distance $\sd$. In particular, the definitions that were previously given for smooth manifolds with a continuous metric apply to this setting.

            \smallskip
            \textbf{Step 1.} We show that, for every $x \in \X \setminus S_2$, there exists $r>0$ such that $\X$ satisfies $\mathrm{CD}(-(n-1)\delta,n)$ in $B_r(x)$. 

            \smallskip
            Suppose first that the point $x \in \X \setminus S_2$ does not belong to $S_1$. In this case, there exists a convex ball $B_r(x)$ not intersecting $S_1 \cup S_2$. 
            Consider two probability measures $\mu_1=\rho_1 \aH^n,\mu_2=\rho_2 \aH^n \in \mathcal{P}_2(B_r(x))$. Any Wasserstein geodesic $\{\xi_t\}_{t \in [0,1]} \subset \mathcal{P}_2(\X)$ between $\mu_1$ and $\mu_2$ is concentrated on geodesics connecting points in $B_r(x)$. These geodesics are themselves contained in $B_r(x)$ by convexity, so that they lie in the region where $\X$ is a smooth manifold with Ricci bounded below by $-(n-1)\delta$. Hence, $\{\xi_t\}_{t \in [0,1]}$ satisfies the required convexity condition of the entropy by \cite{KTSTURM2005}.

            Suppose now that $x \in S_1$.     
            If $R>0$ is small enough,
            using the same arguments of \cite{schlichting2012gluing} (here we use that $\mathrm{I\!I}^M_{S_1}+\mathrm{I\!I}^N_{\partial N} \geq 0$, see also \cite{PerelmanGluingRic}), we obtain the following. 
            There exists a sequence of smooth Riemannian metrics $g_k$ on $\bar{B}_{R}(x) \subset \X$, converging uniformly to $g^\X$ on $\bar{B}_R(x)$, such that $\Ric_{(B_{R}(x),g_k)} \geq -(n-1)\delta -1/k$. 
            We set $r:=R/8$. If $k$ is large enough, geodesics in $(\bar{B}_{R}(x),g_k)$ connecting points in $B_r(x)$ are contained in $B_{R/2}(x)$. In particular, $(\bar{B}_{R}(x),g_k)$ satisfies $\mathrm{CD}(-(n-1)\delta-1/k,n)$ in $B_r(x)$ by the argument used in the previous case. We now consider $(\bar{B}_{R}(x),g^\X)$.

            We denote by $V$ and $V_k$ the volume measures relative to $g^\X$ and $g_k$ in $\bar{B}_R(x)$. 
            Fix two probability measures $\mu_1, \mu_2 \in \mathcal{P}(B_r(x))$ that are absolutely continuous w.r.t. $V$ (or, equivalently, $V_k$) and have compactly supported continuous densities $\mu_1/V$ and $\mu_2/V$. Thanks to \cite[Corollary $29.22$]{Vil}, it is enough to verify the convexity property of the $\mathrm{CD}(-(n-1)\delta,n)$ condition under this extra continuity assumption.
            
            For every $k$, let $\{\xi^k_t\}_{t \in [0,1]}$ be a Wasserstein geodesic between $\mu_1$ and $\mu_2$ in $(\bar{B}_R(x),g_k)$, and let $\Pi_k$ be an optimal plan between the same measures such that, for every $t \in [0,1]$ and every $N \geq n$, it holds
            \[
            U_N(\xi^k_t|V_k) \leq - T^{(t)}_{-(n-1)\delta-1/k,N}(\Pi_k|V_k).
            \]
            By \cite[Theorem $28.9$]{Vil}, the following holds.
            There exists a Wasserstein geodesic $\{\xi_t\}_{t \in [0,1]}$ between $\mu_1$ and $\mu_2$ in $(\bar{B}_R(x),g^\X)$ which arises as limit of $\{\xi^k_t\}_{t \in [0,1]}$. 
            The measures $\Pi_k$ converge weakly to an optimal plan $\Pi$ between $\mu_1$ and $\mu_2$ in $(\bar{B}_R(x),g^\X)$.

            Since $g_k \to g_\X$ uniformly on $\bar{B}_R(x)$, we also have the following.
            The Riemannian distances $\sd_k$ induced by $g_k$ on $\bar{B}_R(x) \times \bar{B}_R(x)$ converge uniformly to the distance $\sd'$ induced by $g^\X$ on $\bar{B}_R(x) \times \bar{B}_R(x)$, which coincides with $\sd$ on $B_r(x) \times B_r(x)$ by the choice of $r$.
            The densities of $V_k$ w.r.t. $V$ converge uniformly to $1$ on $\bar{B}_R(x)$. 
            
            Hence, for every $t \in [0,1]$ and every $N \geq n$, we have uniform convergence
            \begin{align*}
            & \tau^{t}_{-(n-1)\delta-1/k,N}
            \circ \sd_k \to \tau^{t}_{-(n-1)\delta,N}\circ \sd \quad \text{in } B_r(x) \times B_r(x),
            \\
            & (\mu_1/V_k)  \to (\mu_1/V) \quad \text{in } B_r(x),
            \\
            & (\mu_2/V_k)  \to (\mu_2/V) \quad \text{in } B_r(x).
            \end{align*}
            Combining this with the weak convergence of $\Pi_k$ to $\Pi$, we deduce
            \[
            T^{(t)}_{-(n-1)\delta-1/k,N}(\Pi_k|V_k) \to T^{(t)}_{-(n-1)\delta,N}(\Pi|V) \quad \text{for all } t \in [0,1] \text{ and  } N \geq n.
            \]
            Hence, by \cite[Theorem $29.20$]{Vil}, it holds 
            \[
            U_N(\xi_t|V) \leq - T^{(t)}_{-(n-1)\delta,N}(\Pi|V)
            \quad \text{for every } t \in [0,1] \text{ and every } N \geq n.
            \]
            Finally, since geodesics in $(\bar{B}_R(x),g^\X)$ connecting points in $B_r(x)$ are contained in $B_{R/2}(x)$, the Wasserstein geodesic $\{\xi_t\}_{t \in [0,1]}$ in $(\bar{B}_R(x),g^\X)$ is also a  Wasserstein geodesic in $\X$.
            Hence, for every $x \in \X \setminus S_2$, there exists $r>0$ such that $\X$ satisfies $\mathrm{CD}(-(n-1)\delta,n)$ in $B_r(x)$. 
            
            \smallskip
            \textbf{Step 2.} We show that, for every $x \in \X \setminus S_2$, there exists $r>0$ such that geodesics of $\X$ contained in $B_r(x)$ do not branch.

            \smallskip
            If $x \in \X \setminus (S_1 \cup S_2)$, there exists a smooth neighbourhood of $x$, implying the claim. 
            Hence, suppose that $x \in S_1$. 
            Since $M$ and $N$ are smooth manifolds, there exists $k' \in \bb{R}$ such that both their sectional curvatures in a neighbourhood of $x$ are bounded below by $k'$. Hence,
            using the same arguments of \cite{GluingSectionalCurvature} (and using again the condition on the second fundamental forms), we obtain the following. There exist $R>0$ and a sequence of Riemannian metrics $g_k$ on $\bar{B}_{R}(x) \subset \X$, converging uniformly to $g^\X$ on $\bar{B}_R(x)$, such that $\ssf{Sec}_{(B_{R}(x),g_k)} \geq k'-1$.
            
            Hence, for $k$ fixed, every point $y \in B_R(x)$ has a neighbourhood $U_y \subset B_R(x)$, depending on $k$, where triplets of points satisfy the comparison property of Alexandrov spaces with curvature bounded below by $k'-1$ (w.r.t. the distance induced by $g_k$ in $\bar{B}_{R}(x)$).
            This implies, by the proof of Toponogov's Theorem (see \cite[Theorem $10.3.1$]{Coursemetricgeometry}), the existence of $0<r \ll R$ (this time independent of $k$), such that triplets of points in $B_r(x)$ satisfy the comparison property of Alexandrov spaces with curvature bounded below by $k'-1$ (w.r.t. the distance induced by $g_k$ in $\bar{B}_{R}(x)$). 
            
            It is easy to check that since $g_k \to g^\X$ uniformly, then also $(\bar{B}_{R}(x),g^\X)$ satisfies the same triangle comparison property for points contained in $B_r(x)$. This implies that geodesics in $B_r(x)$ do not branch (see \cite{YuBurago_1992}).
            
            \medskip
             \textbf{Step 3.} We first  prove the Laplacian bounds $1$ and $2$ simultaneously, and then we prove $3$.

             \medskip
             Assume that $H_{S_2} \geq \mp \delta$.
            Since $\Ric_{M} \geq -(n-1) \delta$, it holds $$\Delta \sd_{S_2} \leq \pm (n-1) \delta$$ in a neighbourhood of $S_2$ in classical sense (see, for instance, \cite{warpedsplit}). 
            We now use a globalization technique that was previously used in \cite{Weak} and \cite{ketterer2023rigidity}. 
            
            Consider the partition of the transport set of $\sd_{S_2}$ in transport rays $\{r_\alpha\}_{\alpha \in Q}$, 
            $Q$ being a set of indeces with a measure $\ssf{q}$ induced by the partition.
            By \cite[Theorem $3.4$]{CM18}, there exists a disintegration of the measure $\aH^n$ into measures $\{h_\alpha\}_{\alpha \in Q}$, each concentrated on the corresponding transport ray.            
            By \cite[Theorem $4.2$]{MondinoInventiones}
            (the proof works in our setting since the space $\X$ is non-branching and $\mathrm{CD}(-(n-1)\delta,n)$ around every point $x \in \X \setminus S_2$), $\ssf{q}$-almost every measure $h_\alpha$ is absolutely continuous w.r.t. the Lebesgue measure of the corresponding transport ray, it admits a locally Lipschitz continuous density (that we still denote by $h_\alpha$), and it satisfies in weak sense
            \begin{equation} \label{E13}
            (h_\alpha^{\frac{1}{n-1}})''-\delta h_\alpha^{\frac{1}{n-1}} \leq 0.
            \end{equation}
            By \cite[Corollary $4.16$]{CM18} (as before, the proof works in our setting thanks to the previous steps), it holds
            $$
            \Delta \sd_{S_2} = (\log h_\alpha)'+[\Delta \sd_{S_2}]^{sing}
             \quad \text{on } \X \setminus S_2$$ in distributional sense, where $[\Delta \sd_{S_2}]^{sing}$ is a negative singular measure and the transport rays of $\sd_{S_2}$ are now oriented from $S_2$ to $\X \setminus S_2$.
            \par
            Hence, $\ssf{q}$-almost every density $h_\alpha$ satisfies $(\log h_\alpha)' \leq \pm(n-1)\delta$ near $S_2$. Setting $$f_\alpha:=h_\alpha^{\frac{1}{n-1}},$$ we obtain that on the corresponding transport ray $r_\alpha$, thanks to \eqref{E13}, it holds in weak sense
            \[
            \begin{cases}
                f_\alpha''-\delta f_\alpha \leq 0 & \text{on } r_\alpha \\
                f_\alpha' \leq \pm \delta f_\alpha & \text{in a neighbourhood of } S_2 \cap r_\alpha \text{ in } r_\alpha.
            \end{cases}
            \]
            By Riccati Comparison (see, for instance, \cite[Lemma $3.9$]{ketterer2023rigidity}), $f_\alpha$ then satisfies $(\log f_\alpha)' \leq \pm \delta$. In particular, $$\Delta \sd_{S_2} \leq (\log h_\alpha)'=(n-1)(\log f_\alpha)'\leq \pm (n-1)\delta.$$
            We now consider point $3$. 
            As before, we consider the partition of the transport set of $\sd_{S}$ in transport rays $\{r_\alpha\}_{\alpha \in Q}$, and the disintegration of the measure $\aH^n$ into measures $\{h_\alpha\}_{\alpha \in Q}$, each concentrated on the corresponding transport ray.  By our assumptions, if a point $p$ of a transport ray $r_\alpha$ 
            lies in the $D$-neighbourhood of $S$ in $\X \setminus N_S$, then all the points of $r_\alpha$ between $p$ and $S$ lie in this region.
            Hence, the same argument that we used for points $1$ and $2$ can be repeated.
            \end{proof}
\end{lemma}

The next lemma is standard and it follows from \cite{warpedsplit} and \cite{Kasue}. We sketch a proof for the sake of completeness.

\begin{lemma}\label{CL2}
    Let $\delta \in \{0,1\}$ and let $(M^n,g)$ be a Riemannian manifold with $\Ric_M \geq -(n-1)\delta$. Let $S_1$ be a union of connected components of $\partial M$ and let $S_2:=\partial M \setminus S_1 \neq \emptyset$. Let $H_{S_1} \geq \delta$, $H_{S_2} \geq -\delta$ and let $\sd_{S_1}+\sd_{S_2}$ be constant. 
    Then, there exists $l>0$ such that $M$ is isometric to $S_1 \times [0,l]$ with
    the warped product metric $\de s^2 = e^{-2\delta t}g_1 + \de t^2$, where $g_1$ is the metric on $S_1$.
    \begin{proof}
        Since $\sd_{S_1}+\sd_{S_2} \equiv D $ is constant, every point of $M$ lies on a unique geodesic realizing the distance from $S_1$ to $S_2$. In particular, $\sd_{S_1}$ and $\sd_{S_2}$ are both smooth and they satisfy $$\Delta \sd_{S_1} \leq -(n-1)\delta \quad \text{and}\quad  \Delta \sd_{S_2} \leq (n-1)\delta.$$ Hence, $\Delta \sd_{S_1} = -(n-1)\delta$.
        Using $t=\sd_{S_1}$ as a coordinate on $M$, we can write $g=\de t^2+g^\perp(t)$, where $g^\perp(t)$ is a metric on the $t$ enlargement of $S_1$ in $M$.
        
        Using Cauchy-Schwarz inequality first, and then plugging $\sd_{S_1}$ into the Bochner Formula, we obtain
        \[
        (n-1)\delta=\frac{(\Delta \sd_{S_1})^2}{n-1} \leq |\Hess(\sd_{S_1})|^2=-\Ric_M(\nabla \sd_{S_1},\nabla \sd_{S_1}) \leq (n-1)\delta.
        \]
        In particular, all inequalities in the previous expression are equalities, forcing $\Hess(\sd_{S_1})=-\delta g^{\perp}$ everywhere on $M$. Denoting by $L_{\partial_t}$ the Lie derivative in the direction of $\nabla \sd_{S_1}$, the previous equality implies 
        \begin{equation} \label{E12}
        L_{\partial_t} \, g^\perp=-2\delta g^\perp.
        \end{equation}
        Consider now the map $\phi:S_1 \times (0,D) \to M \setminus (S_1 \cup S_2)$ which sends $(x,t)$ to the point at distance $t$ from $x$ along the flow line of $\nabla \sd_{S_1}$. This map is bijective since every point of $M$ lies on a unique geodesic realizing the distance from $S_1$ to $S_2$. Equipping $S_1 \times (0,D)$ with the metric $\de s^2=\de t^2+e^{-2\delta t}g_1$, we obtain that $\phi$ is an isometry by \eqref{E12} and the definition of Lie derivative, concluding the proof.
    \end{proof}
\end{lemma}

\bigskip

The next lemmas deal with \emph{sets minimizing the perimeter} in Riemannian manifolds, and are needed to prove Theorem \ref{CT4}. In Euclidean space, the subject is by now classical, see for instance the monographs \cite{Giusti, Maggi}. When the ambient space is not $\bb{R}^n$, many properties have been obtained in the more general setting of metric measure spaces. An account of this theory can be found in \cite{Weak}.
\medskip

Let $(M,g)$ be a Riemannian manifold, and let $E \subset M$ be a Borel set. Given an open set $A \subset M$, the perimeter of $E$ relative to $A$ is defined as
$$
P(E, A) := \inf \left\{ \liminf_{k \to \infty} \int_A |\nabla f_k|\, \de \ssf{Vol} \ : \ f_k \in C^\infty(A)
, \ f_k \to \chi_E\ \mbox{ in } L^1_\mathrm{loc}(A)\right\}.
$$
The set $E \subset M$ is said to have locally finite perimeter if $P(E, B_r(x)) < + \infty$ for all $x \in M$ and $r > 0$.
If $E$ has locally finite perimeter, there exists a unique Radon measure $\mu$ such that $\mu(A)=P(E,A)$ if $A \subset M$ is open. This measure is denoted $P(E,\cdot)$.
\medskip

A set of locally finite perimeter $E \subset M$ is \emph{perimeter minimizing} in an open set $A \subset M$ if, for every $F \subset M$ such that $F \Delta E \subset \subset A$, it holds that $P(E,A) \leq P(F,A)$. The set $E \subset M$ is \emph{locally perimeter minimizing} in $A$ if for every $x \in A$ there exists $r_x>0$ such that $E$ is perimeter minimizing in $B_{r_x}(x)$.
\medskip

If $E \subset M^n$ locally minimizes the perimeter in $M$, then it admits both an open and a closed representative, and these have the same topological boundary, which is a smooth hypersurface outside of a set $\Sigma \subset \partial E$ of Hausdorff dimension at most $n-8$. Whenever we refer to the boundary of a locally perimeter minimizing set, we mean the topological boundary of its open (or closed) representative.
\medskip 

Lemma \ref{L|Morgan} concerns a minimizing property of the boundaries of mean-convex sets in Riemannian manifolds. It is well known that minimal hypersurfaces locally minimize area (see \cite{LG96}). Lemma \ref{L|Morgan} can be viewed as an analogue of this fact for sets with mean-convex boundary, and its proof follows by adapting the argument in \cite{LG96}.  Although this result is likely known to experts, we include the proof for completeness.

\begin{lemma} \label{L|Morgan}
    Let $(M^n,g)$ be a smooth Riemannian manifold, let $E \subset M$ be a smooth open set with mean-convex boundary, and let $x \in \partial E$. Then there exists a ball $B$ centered in $x$ such that, for every $E \subset E' \subset M$ with $E' \Delta E \subset \subset B$, it holds that $P(E,B) \leq P(E',B)$.
    \begin{proof}
        Let $B$ be a small ball of $M$ centered in $x$. Let $C \subset \partial E$ be an $(n-2)$-dimensional submanifold of $\partial E$ containing $x$, and let $P_0:\partial E \cap B \to C$ be the closest point projection inside $M$ restricted to $\partial E$. Observe that if $B$ is small enough, this map is well defined and smooth.
        Using the smooth coarea formula (see \cite[Exercise III.12 point c]{Chavel}) w.r.t.\ $P_0:\partial E \cap B \to C$, we can rewrite the area of $\partial E \cap B$ as
        \begin{equation} \label{E1LawlorMorgan}
        \aH^{n-1} (\partial E \cap B)=\int_C \int_{B \cap \partial E_c} 1/J_{n-2}(P_{0}) \,\de \aH^{1} \, \de \aH^{n-2}(c), 
        \end{equation}
        where $\partial E_c:=P_0^{-1}(c) \subset \partial E \cap B$ and $J_{n-2}(P_0)$ is the appropriate coarea factor, 
        i.e.\ for every $x \in \partial E \cap B$
        \[
        J_{n-2}(P_0)(x):=|\mathrm{det}({d_xP_0}_{|\mathrm{Ker}(d_xP_0)^{\perp,T_x\partial E}})|,
        \]
        where 
        $\mathrm{Ker}(d_xP_0)^{\perp,T_x\partial E}$ denotes the orthogonal complement of $\mathrm{Ker}(d_xP_0)$ in $T_x \partial E$.
        Let now $P_1: B \to \partial E$ be the natural nearest point projection in $M$. Modulo shrinking $B$, also $P_1$ is well defined and smooth. Set now $P:=P_0 \circ P_1:B \to C$ and, given a point $y \in \partial E \cap B$, consider the coarea factor $J_{n-2}(P)(y)$.
        We next relate $J_{n-2}(P)(y)$ to $J_{n-2}(P_0)(y)$.
        Observe that, using the identification $T_y \partial E \subset T_y M$, it holds 
        \[
        \mathrm{Ker}(d_y P)=\mathrm{Ker}(d_y P_0) \oplus T_y\partial E^{\perp, T_yM},
        \]
        so that        
        $\mathrm{Ker}(d_y P)^{\perp,T_yM}=\mathrm{Ker}(d_y P_0)^{\perp, T_y \partial E}$.
        Moreover, by the chain rule and Binet's Formula, it holds
        \begin{align*}
        J_{n-2}(P)(y)& :=
        |\mathrm{det}({d_yP}_{|\mathrm{Ker}(d_yP)^{\perp,T_yM}})|
        \\
        & =
        |\mathrm{det}({d_yP_1}_{|\mathrm{Ker}(d_y P_0)^{\perp, T_y \partial E}})|
        |\mathrm{det}({d_yP_0}_{|\mathrm{Ker}(d_y P_0)^{\perp, T_y \partial E}})| \\
        & =J_{n-2}(P_0)(y),
        \end{align*}
        where the last identity follows from the fact that $\mathrm{Ker}(d_y P_0)^{\perp, T_y \partial E} \subset T_y\partial E$.
        Hence, we can rewrite \eqref{E1LawlorMorgan} as
        \begin{align} \label{E|Morgan1}
        \aH^{n-1} (\partial E \cap B)
        & =
        \int_C \int_{B \cap \partial E_c} 1/J_{n-2}(P) \,\de \aH^{1} \, \de \aH^{n-2}(c)
        \\ \nonumber
        & =
        \int_C \mathrm{length}_c(\partial E_c \cap B) \,\de \aH^{n-2}(c),
        \end{align}
        where $\mathrm{length}_c$ denotes the length of a curve in the smooth surface $\{P=c\} \subset B$ with conformally rescaled metric $J_{n-2}(P)^{-2}g$.

        Observe now that since every point in the image of $P$ is a regular value, then $\partial E_c$ is a smooth curve, and it is equal to the boundary of $$E_c:=E \cap \{P=c\},$$ inside the smooth surface $\{P=c\}$. 
        \smallskip
        
        \textbf{Key Step:}  $E_c$ is convex in $\{P=c\}$, endowed with the metric $J_{n-2}(P)^{-2}g$. 
\\        Assume by contradiction that this is not the case, so that there exists $c_0 \in C$ such that $\partial E_{c_0}$ has a point of strictly negative (mean) curvature in $\{P=c_0\}$ with the metric $J_{n-2}(P)^{-2}g$. By continuity, there is a neighbourhood $U \subset B$ of $x_0$ in $M$, such that every point of $\partial E_c \cap U$ for some $c \in C$ has strictly negative mean curvature in the corresponding rescaled metric of $\{P=c\}$. Let now $X \in TU \subset TM$ be a smooth vector field such that: 
        \begin{enumerate}
            \item $X_{|\partial E \cap U}$ is non-vanishing and outer pointing;
            \item if $p \in \partial E_c \cap U$ for some $c \in C$, then $X(p) \in T_p \{P=c\} \subset T_pM$.
        \end{enumerate}
        It is clear that, modulo restricting $U$, a vector field with these properties exists. Consider now a nonnegative function $\phi \in C^{\infty}_c(\partial E \cap U)$ which is equal to $1$ in a neighbourhood of $x_0$, and consider for every $t>0$ small enough the flow $\phi_t(\cdot):\partial E \to M$ w.r.t. $\phi X$. Set $\Sigma^t:=\phi_t(\partial E)$, and observe that these only differ from $\partial E$ in a compact set of $U$, and converge smoothly to $\partial E$. Moreover, setting $\Sigma_c^t:=\Sigma^t \cap \{P=c\}$, one has that $\phi_t(\partial E_c)=\Sigma^t_c$ by the properties of $X$.

        By the first variation of the area w.r.t.\ compactly supported variations (see \cite[Chapter 1]{LiG}), we have
        \begin{equation}
            \frac{d}{dt} \aH^{n-1}(\Sigma^t \cap U)_{|t=0} =\int_{\partial E} \phi H_{\partial E} \cdot X \, d\aH^{n-1} \geq 0,
        \end{equation}
        where $H_{\partial E}$ is the mean curvature vector, i.e.\ a normal vector to $\partial E$ whose modulus is the modulus of the mean curvature of $\partial E$ in $E$ (with the understanding that this vector points outwards if the mean curvature is non-negative, and inwards otherwise).

        At the same time, using the coarea formula as in the first part of the proof, and denoting by $\mathrm{length}^t_c(\Sigma^t_c \cap U)$ the length of the curve
        $\Sigma^t_c \cap U$
        with respect to the conformally rescaled metric $J^{-2}_{n-2}({P}_{|\Sigma^t})g$ on $\Sigma^t$,
        one has
        \begin{align*}
             \aH^{n-1}(\Sigma^t \cap U) &=
              \int_C \mathrm{length}^t_c(\Sigma^t_c \cap U) \,\de \aH^{n-2}(c)
              \\
              & =
              \int_C \mathrm{length}^t_c(\phi_t(\partial E_c \cap U)) \,\de \aH^{n-2}(c).
        \end{align*}
        We now claim that
        \begin{equation} \label{E|lawsonMorgan2}
            \frac{d}{dt}\mathrm{length}^t_c(\phi_t(\partial E_c \cap U))_{|t=0}=
            \frac{d}{dt}\mathrm{length}_c(\phi_t(\partial E_c \cap U))_{|t=0}.
        \end{equation}
        Assume for the moment that the claim is true.
        In this case, passing the derivatives under the integral sign, using the first variation of the area in each slice ${P=c}$, and denoting by $H^c_{\partial E_c}$ the mean curvature of $\partial E_c$ in $\{P=c\}$ with the conformally rescaled metric, we get
        \begin{align*}
            \frac{d}{dt} \aH^{n-1}(\Sigma^t \cap U)_{|t=0}& =\int_C
            \frac{d}{dt}\mathrm{length}_c(\phi_t(\partial E_c \cap U))_{|t=0} \, d\aH^{n-2} (c) \\
            &=  \int_C \int_{\partial E_c} \phi  (H^{c}_{\partial E_c} \cdot X) \, {J^{-1}_{n-2}}(P)d\aH^1 d\aH^{n-2} (c)<0.
        \end{align*}
        The final strict inequality comes from our mean curvature assumptions on $\partial E_c \cap U$ inside $\{P=c\}$ with the rescaled metric, giving the desired contradiction.
        
        To conclude the proof of the key step, we need to prove \eqref{E|lawsonMorgan2}.
         Observe that as $t \to 0$, we have that $J_{n-2}(P_{|\Sigma^t}) \to J_{n-2}({P}_{|\partial E})$ smoothly (in the natural sense, pulling back everything to $\partial E$ via $\phi_t$), and that $J_{n-2}({P}_{|\partial E})$ is strictly positive.
         Hence, it suffices to show that
        \begin{equation}
            \sup_{U \cap \Sigma^t}\frac{|{J_{n-2}}(P_{|\Sigma^t})-J_{n-2}(P)|}{t} \to 0.
        \end{equation}
        For every $x \in \Sigma^t_c$, since $\mathrm{Ker}(d_xP)=T_x\{P=c\}$,
        we have
        \begin{equation}
            \mathrm{Ker}(d_xP) \cap T_x \Sigma^t=T_x \Sigma^t_c,
        \end{equation}
        so that
        \begin{equation}
           (\mathrm{Ker}(d_xP) \cap T_x \Sigma^t)^{\perp,T_x \Sigma^t}=T_x{\Sigma_c^t}^{\perp,T \Sigma^t}.
        \end{equation}
        Denote by $\Pi_{t,x}:T_x{\Sigma_c^t}^{\perp,T \Sigma^t} \to \mathrm{Ker}(d_xP)^{\perp,T_x M}$ the restriction of the orthogonal projection inside $T_xM$.
        It holds,
        \begin{align*}
            J_{n-2}(P_{|\Sigma^t})(x)
            & 
            =
            |\mathrm{det}({d_xP}_{|(\mathrm{Ker}(d_xP) \cap T_x\Sigma^t)^{\perp,T_x \Sigma^t}})| \\
            & = 
            |\mathrm{det}({d_xP}_{|\mathrm{Ker}(d_xP)^{\perp,T_x M}} \circ \Pi_{t,x})|
             \\
            & = 
            J_{n-2}(P)(x)
            |\mathrm{det} (\Pi_{t,x})|.
        \end{align*}
       Hence, it suffices to show that
       \begin{equation} \label{E|lawsonMorgan3}
           \sup_{x \in \Sigma^t \cap U} \frac{\big| |\mathrm{det} (\Pi_{t,x})|-1 \big|}{t} \to 0 \quad \text{as } t \to 0.
       \end{equation}
       This follows because for every $y \in \partial E \cap U$, setting $y_t:=\phi_t(y)$, the function $t \to \mathrm{det} (\Pi_{t,y_t})$ is smooth in $[0,\epsilon]$ for some $\epsilon>0$, and it is a product of the cosines of the angles formed by (an orthonormal basis of) $T_{y_t}{\Sigma_c^t}^{\perp,T \Sigma^t}$ with $\mathrm{Ker}(d_{y_t}P)^{\perp,T_{y_t} M}$. These angles vary smoothly, and they are all equal to $0$ at $t=0$. Expanding the cosine, we obtain that $\mathrm{det} (\Pi_{t,y_t})$ has vanishing first derivative in $0$, which then implies 
       \begin{equation} \label{E|lawsonMorgan5}
            \frac{\big| |\mathrm{det} (\Pi_{t,y_t})|-1 \big|}{t} \to 0  \quad \text{for every } y \in \partial E \cap U \text{ as } t \to 0.
       \end{equation}
        Moreover, there is a compact set $K \subset \partial E \cap U$ such that $y_t=\phi_t(y)=y$ for every $y \in \partial E \cap U \setminus K$, so that \eqref{E|lawsonMorgan3} amounts to showing that
        \begin{equation}
            \label{E|lawsonMorgan4}
           \sup_{y \in K} \frac{\big| |\mathrm{det} (\Pi_{t,y_t})|-1 \big|}{t} \to 0  \quad \text{as } t \to 0.
        \end{equation}
        Observe now that the function $(y,t) \mapsto (\mathrm{det} (\Pi_{t,y_t})-1)/t$  defined on $(\partial E \cap U) \times [0,\epsilon]$ for some $\epsilon>0$ is also smooth (because $(y,t) \mapsto \mathrm{det} (\Pi_{t,y_t})$ is smooth, and it is equal to $1$ in $t=0$). Then, \eqref{E|lawsonMorgan4} follows from \eqref{E|lawsonMorgan5} via a compactness argument. This concludes the proof of the key step.
\smallskip

        We now complete the proof of the lemma. By the Key Step, each set $E_c \subset \{P=c\}$ has convex boundary. By a standard property of convex sets (which follows by combining \cite{BishopConvexity} and \cite[Lemma 1 Section 3]{Alexander}), for every $c \in C$ and $x \in \partial E_c$, there exists a neighbourhood $U^c_x$ of $x$ in $\{P=c\}$,
        such that any curve of $\{P=c\}$ which lies outside of $E_c$ and connects two points of $\partial E_c \cap U^c_x$ is longer than the corresponding portion of $\partial E_c$ with the same endpoints.
        The neighbourhoods $U^c_x \subset \{P=c\}$ vary continuously as $c$ and $x$ vary. Hence, having fixed $x \in \partial E$, there exists a neighbourhood $U_x$ in $M$ such that for every $c \in C$ and $y \in \partial E_c \cap U_x$ one can choose $U^c_y:=U_x \cap \{P=c\}$. 
        
        Fix now a ball $B \subset U_x$ containing $x$ and 
        let $E'\subset M$ be a smooth set such that $E \subset E'$ and $E' \Delta E \subset \subset B$.
        Denoting $\partial E'_c:=E' \cap \{P=c\}$, it holds
        \begin{align*} 
        \aH^{n-1} (\partial E' \cap B)
        =
        \int_C \int_{\partial E_c' \cap B}J_{n-2}^{-1}(P_{|\partial E'}) \, d\aH^1 d\aH^{n-2}(c) \\
         \geq
        \int_C \mathrm{length}_c(\partial E'_c \cap B) \,d \aH^{n-2}(c),
        \end{align*}
        where the last inequality is justified by the fact that $J_{n-2}(P_{|\partial E'}) \leq J_{n-2}(P)$ immediately from the definition. Using that $\mathrm{length}_c(\partial E'_c \cap B) \geq \mathrm{length}_c(\partial E_c \cap B)$ and \eqref{E|Morgan1}, we then deduce
        \begin{align*} 
        \aH^{n-1} (\partial E' \cap B)
        & \geq
        \int_C \mathrm{length}_c(\partial E'_c \cap B) \,d \aH^{n-2}(c)
        \\
        & \geq \int_C \mathrm{length}_c(\partial E_c \cap B) \,d \aH^{n-2}(c)
        = \aH^{n-1} (\partial E \cap B).
        \end{align*}
        In the case when $E' \subset M$ is not smooth, the statement follows by an approximation argument combined with the previous part of the proof.
    \end{proof}
\end{lemma}

The classical
 \cite[lemma $4$, section $11$]{MeeksSimonYau} states that a manifold with mean convex disconnected boundary and dimension $ \leq 7$ contains a properly embedded stable minimal hypersurface. The next lemma is a variant of this result. We consider a manifold $M$ with mean convex disconnected boundary and arbitrary dimension. We show that either one connected component of the boundary is a stable minimal hypersurface, or $M \setminus \partial M$ contains a non-trivial perimeter minimizing set.

 The logic behind this duality is that if the minimal hypersurface given by Meeks-Schoen-Yau touches the boundary, then it coincides with a boundary component by the maximum principle. Otherwise, since it is constructed via a minimization procedure, 
 it will bound a set which will minimize the perimeter. Lemma \ref{L|yau} also gives some area bounds for the constructed perimeter minimizing set. These are only needed for the conclusive remarks in Section \ref{S5}.

\begin{lemma} \label{L|yau}
    Let $(M^n,g)$ be a Riemannian manifold with mean-convex disconnected boundary $\partial M$, and let $\Sigma_1 \subset \partial M$ be a connected component. Then, at least one of the following two assertions holds.
    \begin{enumerate}
        \item \label{item1Yau} One connected component $\Sigma \subset \partial M$ of $\partial M$ is minimal,  stable, and it satisfies 
        \[
        \aH^{n-1}(\Sigma \cap \bar{B}_r(p)) \leq \aH^{n-1}(\partial B_r(p))+\aH^{n-1}(\bar{B}_r(p) \cap \Sigma_1)
        \]
        for every $p \in M$ and every $r >0$.
        \item \label{item2Yau}
        There exists an open set $E \subset M$ such that:
        \begin{itemize}
            \item $E$ minimizes the perimeter in $M \setminus \partial M$.
            \item $\partial E \neq \emptyset$, $\partial E \cap \partial M=\emptyset$, and $E \supset \Sigma_1$.
            \item $\aH^{n-1}(\partial E \cap \bar{B}_r(p)) \leq \aH^{n-1}(\partial B_r(p))+\aH^{n-1}(\bar{B}_r(p) \cap \Sigma_1)$ for every $p \in M$ and every $r >0$.
        \end{itemize}
    \end{enumerate}
    \begin{proof}
        Set $\Sigma_2 :=\partial M \setminus \Sigma_1$. Let $(\tilde{M},\tilde{g})$ be an enlargement of $M$ beyond its boundary i.e.,\ an open manifold containing open sets $E_1,E_2 \subset \tilde{M}$ with $\partial E_1=\Sigma_1$, $\partial E_2=\Sigma_2$, and $\tilde{M} \setminus (E_1 \cup E_2)=M$ isometrically. Unless otherwise specified, in the rest of the proof, all balls and boundaries are taken in $\tilde{M}$.
        Let $p \in M$, and for $i \in \bb{N}$ large enough consider the minimization problem
        \begin{equation} \label{E|var1}
        \min \{P(E,B_i(p)):E \supset E_1, \, E \cap E_2=\emptyset, \, E=E_1 \text{ in }B_i(p) \setminus B_{i-1}(p)\}.
        \end{equation}
        For every $i \in \bb{N}$ sufficiently large, we denote by $E^i \subset \tilde{M}$ the set realizing the minimum in the previous minimization problem. Such set exists by the direct method of calculus of variations. We construct the set $E \subset M$ in several steps. We then show that either $E$ satisfies item \ref{item2Yau}, or item \ref{item1Yau} holds.

        \textbf{Step 1:} The sequence $E^i$ can be taken to be increasing. 
        
        Assume that the partial sequence $\{E^i\}_{i=1}^l$ is increasing. We need to show that $E^{l+1}$ can be chosen so that $E^{l+1} \supset E^l$. Indeed, let $\tilde{E}^{l+1}$ be a minimum of the variational problem \eqref{E|var1}, and set $E^{l+1}:=\tilde{E}^{l+1} \cup E^l$. The set $E^{l+1}$ is trivially a competitor for the variational problem \eqref{E|var1}. Moreover,
        \begin{align} \label{E|var3}
        \nonumber
        P &(E^{l+1},B_{l+1}(p)) =P(\tilde{E}^{l+1}\cup E^l,B_{l+1}(p)) \\
        &\leq P(\tilde{E}^{l+1},B_{l+1}(p))+P(E^{l},B_{l+1}(p))-P(\tilde{E}^{l+1}\cap E^l,B_{l+1}(p)).
        \end{align}
        In addition, $P(E^{l},B_{l}(p)) \leq P(\tilde{E}^{l+1}\cap E^l,B_{l}(p))$, since $\tilde{E}^{l+1}\cap E^l$ is a competitor for $E^l$ in \eqref{E|var1}. Since $E^l=\tilde{E}^{l+1}\cap E^l$ in $B_{l+1} \setminus B_{l-1}(p)$, we deduce $P(E^{l},B_{l+1}(p)) \leq P(\tilde{E}^{l+1}\cap E^l,B_{l+1}(p))$. Combining with \eqref{E|var3}, we obtain
        \[
        P(E^{l+1},B_{l+1}(p)) \leq P(\tilde{E}^{l+1},B_{l+1}(p)).
        \]
        This concludes the proof of Step 1.

        \textbf{Step 2:} There exists a locally finite covering with balls $\{B_\alpha\}_{\alpha}$ of $\bar{M} \subset \tilde{M}$ such that,  for every $\alpha$, if $i$ is large enough (depending on $\alpha$), the sets $E^i \subset \tilde{M}$ minimize the perimeter in the ball $B_\alpha$.

        It is sufficient to show that for every $x \in \bar{M} \subset \tilde{M}$ there exists a ball $B$ of $\tilde{M}$ containing $x$ such that the sets $E^i$ are perimeter minimizing in $B$. If $x \in M \setminus \partial M$, one can take any ball centered in $x$ which does not intersect $\partial M$. 
        So, we consider the case $x \in \partial M$. We assume that $x \in \Sigma_2$, the case $x \in \Sigma_1$ being analogous.
        Since $\partial M$ is mean-convex, by Lemma \ref{L|Morgan}, there exists a ball $B \subset \tilde{M}$ centered in $x$ such that
        for every $M \subset C \subset \tilde{M}$ with $C \Delta M \subset \subset B$, it holds $P(M,B) \leq P(C,B)$. This ball can be taken small so that it is disjoint from $\Sigma_1$.
        We show that each $E^i$ is perimeter minimizing in the ball $B \subset \tilde{M}$. Let $C \subset \tilde{M}$ be such that $C \Delta E^i \subset \subset B$. Since $E^i$ minimizes the variational problem \eqref{E|var1}, it holds
        \begin{equation} \label{E|var2}
        P(E^i,B) \leq P(C \cap M,B).
        \end{equation}
        By our choice of $B$, it holds $P(M,B) \leq P(C \cup M,B)$. Combining with \eqref{E|var2}, we deduce
        \[
        P(E^i,B) \leq P(C \cap M,B) \leq P(C,B)+P(M,B)-P(C \cup M,B) \leq P(C,B),
        \]
        proving our claim.

        Hence, up to passing to a subsequence, there exists $E \subset \tilde{M}$ locally minimizing the perimeter such that $E^i \to E$ in $\sL^1_{loc}(\tilde{M})$ (see \cite[Theorem 1.19 and Lemma 9.1]{Giusti}). Since each $E^i$ minimizes the perimeter in $B_{i-1}(p) \cap (M \setminus \partial M)$ by construction, $E$ also minimizes the perimeter in $M \setminus \partial M$.
        By construction, it also holds $E \supset E_1$ and $E \cap E_2=\emptyset$.

        \textbf{Step 3:} Let $U \subset \tilde{M}$ be open and bounded. If $C \subset \tilde{M}$ is such that $E \supset C \supset E_1$, and $C \Delta E \subset \subset U$, then $P(E,U) \leq P(C,U)$.

        The claim follows combining Step 1 and a standard argument (see, for instance, \cite[Proposition 1.3]{DephilippisSimonsCone}). We report the argument for the sake of completeness. Consider $C \subset \tilde{M}$ as in the statement of step 3. The sets $C \cap E^i$ satisfy $P(E^i,U) \leq P(C \cap E^i,U)$ if $i$ is large enough. Hence,
        \[
        P(E^i \cup C,U) \leq P(E^i,U)+P(C,U)-P(C \cap E^i,U) \leq P(C,U).
        \]
        The sets $C \cup E^i$ converge in $\sL^1(U)$ to $E$ by Step 1. Hence, by lower semicontinuity, it holds
        $
        P(E,U) \leq P(C,U)
        $, concluding the proof of Step 3.
         
        \textbf{Step 4:} Denoting by $B^M$ and $\partial^M$ respectively balls and boundaries in $M$ (and not in $\tilde{M}$), for every $r>0$, it holds
        \[
        \aH^{n-1}(\partial E \cap \bar{B}^M_r(p)) \leq \aH^{n-1}(\partial^M B^M_r(p))+\aH^{n-1}(\bar{B}^M_r(p) \cap \Sigma_1). 
        \]
        
         Let $U \subset \tilde{M}$ be an open set of $\tilde{M}$ containing $\bar{B}^M_r(p)$.
        By Step 3, it holds 
        \begin{align*}
        \aH^{n-1} (\partial E \cap U) =P(E,U) 
        \leq P((E \setminus B^M_r(p)) \cup E_1,U) \\
         \leq  \aH^{n-1}(\partial [(E \setminus B^M_r(p)) \cup E_1] \cap U).
        \end{align*}
        Taking an infimum over all possible open sets $U \subset \tilde{M}$ such that $U \supset \bar{B}^M_r(p)$, we have
        \[
        \aH^{n-1} (\partial E \cap \bar{B}^M_r(p)) 
         \leq  \aH^{n-1}(\partial [(E \setminus B^M_r(p)) \cup E_1] \cap \bar{B}^M_r(p)).
        \]
        Since
        \[
        \partial [(E \setminus B^M_r(p)) \cup E_1] \cap \bar{B}^M_r(p) \subset \partial^M B^M_r(p) \cup  \Sigma_1,
        \] 
        Step 4 follows.

        We now conclude the proof of the lemma. Assume that $x \in \partial E \cap \partial M$. We wish to show that in this case item \ref{item1Yau} is satisfied. Taking a blow-up of $E$ in $x$, we find a perimeter minimizing set $E_\infty \subset \bb{R}^n$ which does not intersect the half-space of $\bb{R}^n$ whose boundary is given by the blow-up of $\partial M$ in $x$. Hence, $E_\infty \subset \bb{R}^n$ is a half-space (see \cite[Theorem 15.5]{Giusti}). By the regularity theory, it follows that $E$ is a smooth set in a neighbourhood of $x$ in $\tilde{M}$. 
        Since $\partial M$ is mean-convex and $\partial E$ locally lies on the interior side of $\partial M$, by the maximum principle, we infer that $\partial E =\partial M$ in a neighbourhood of $x$. 
        By a standard argument, it follows that $\partial E$ contains a connected component  $\Sigma \subset \partial M$ of $\partial M$. Since $E$ is locally perimeter minimizing in $\tilde{M}$, such connected component $\Sigma \subset \partial E$ has vanishing mean curvature. By Step 3, it follows that $\Sigma$ is a stable minimal hypersurface, and by Step 4 it follows that it satisfies the required volume growth.

        Assume now that $\partial E \cap \partial M=\emptyset$. We wish to show that $E$ satisfies the conditions in item \ref{item2Yau}. Indeed, $E$ is locally perimeter minimizing, so that it has an open representative. By construction, $E$ minimizes the perimeter in $M \setminus \partial M$. Since $E \supset E_1$ and $E \cap E_2=\emptyset$, the boundary $\partial E$ in $\tilde{M}$ is nonempty, and since $\partial E \cap \partial M=\emptyset$ by assumption, we have $\partial^M E \neq \emptyset$ as well. Finally, the desired volume growth condition follows by Step 4.
    \end{proof}
\end{lemma}

\section{Proofs of Theorems \ref{CT1}, \ref{CT1*}, and \ref{CT2}}
    We now prove Theorems \ref{CT1}, \ref{CT1*}, and \ref{CT2}. We recall that a manifold $(M,g)$ without boundary is said to be \emph{parabolic} if it has no positive fundamental solution for the Laplacian. This is equivalent to asking that for every $x \in M$ and $r>0$, it holds
    \[
    \inf \Big\{\int_M |\nabla \phi|^2 \, \de \ssf{Vol}: \phi \in \Lip_c(M): \phi \equiv 1 \text{ on } B_r(x) \Big\}=0.
    \]
    If $\Ric_{M} \geq 0$, parabolicity is equivalent to requiring that $$\int_1^{+\infty} \frac{t}{\ssf{Vol}(B_t(x))} \, \de t=+\infty,$$ for some $x \in M$. For the proofs of these facts, we refer to \cite{GP} and the monograph \cite{LiG}. The following theorem corresponds to Theorems \ref{CT1} and \ref{CT1*}.
    
    \medskip
    \begin{thm} \label{CT1body}
    Let $\delta \in \{0,1\}$ and let $(M^n,g)$ be a Riemannian manifold, with $\Ric_{M} \geq -(n-1) \delta$. Let $S_1$ be a non-empty connected component of $\partial M$, and let $S_2:=\partial M \setminus S_1$ be non-empty as well. Assume that $S_1$ is parabolic, with  $\mathrm{I\!I}_{S_1} \geq \delta$ and $\Ric_{S_1} \geq 0$, while $H_{S_2} \geq -\delta$.
    Then, there exists $l>0$ such that $M$ is isometric to $S_1 \times [0,l]$ with
    the warped product metric $\de s^2 = e^{-2\delta t}g_1 + \de t^2$, where $g_1$ is the metric on $S_1$.
\begin{proof}
    We first introduce some tools.
        Consider the manifold $S_1 \times (0,+\infty)$ equipped with the warped product metric $\tilde{g} = e^{2\delta t}g_1 + \de t^2$. We denote by $\ssf{Vol}$ the Riemannian volume in $S_1 \times (0,+\infty)$ w.r.t. the metric $\tilde{g}$ and by $\ssf{Vol}_{S_1}$ the Riemannian volume on $S_1$ w.r.t. the metric $g_1$.
        
        Let $\sd_{S_1 \times\{0\}}$ be the distance from $S_1 \times \{0\}$ in $S_1 \times (0,+\infty)$ and let $\sd_{S_1}$ be the distance from $S_1$ in $M$.
        By standard computations (see \cite{Pet}), it holds
        $$
        \Delta \sd_{S_1 \times \{0\}}=(n-1)\delta, \quad \text{on } S_1 \times (0,+\infty).$$ 
        Similarly, the second fundamental form of $S_1 \times \{0\}$ in $S_1 \times (0,+\infty)$ satisfies
        $
        \mathrm{I\!I}_{S_1 \times \{0\}}=-\delta
        $. Moreover, $\Ric_{S_1} \geq 0$ together with Gauss' equation implies that $\Ric_{S_1 \times (0,+\infty)} \geq -(n-1)\delta$. 
        \par 
        Consider the metric space obtained by glueing $M$ and $S_1 \times (0,+\infty)$ along the two boundaries isometric to $S_1$. 
        On this metric space, we consider the distance function from $S_2$ and we denote it by $\sd_{S_2}$. 
        In particular, this allows us to define $\sd_{S_2}$ on $S_1 \times (0,+\infty)$ (seen as a subset of the previously defined metric space). 
        We claim that the function $\bar{\sd}:=\sd_{S_2}-\sd_{S_1 \times \{0\}}$ is constant on $S_1 \times (0,+\infty)$.
        
    Let us first show how this is enough to conclude the proof. If $\bar{\sd}:=\sd_{S_2}-\sd_{S_1 \times \{0\}}$ is constant on $S_1 \times (0,+\infty)$, then $\sd_{S_2}$ is constant when restricted to ${S_1}$, which in turn implies that there exists a geodesic $\gamma \subset M$ joining $S_1$ and $S_2$ whose length is equal to the distance between $S_1$ and $S_2$. By triangle inequality, $\sd_{S_1}+\sd_{S_2}:M \to \bb{R}$ attains its minimum on any point of $\gamma$, so that, being superharmonic by Lemma \ref{l|Laplacian in the glued space}, it is constant by the Maximum Principle. The conclusion then follows by Lemma \ref{CL2}. 

    Hence, we only need to prove that $\bar{\sd}:=\sd_{S_2}-\sd_{S_1 \times \{0\}}$ is constant on $S_1 \times (0,+\infty)$. 
        
        By Lemma \ref{l|Laplacian in the glued space}, it holds $\Delta \sd_{S_2} \leq (n-1)\delta$ in distributional sense on $S_1 \times (0,+\infty)$. Hence, 
        \begin{equation} \label{E1}
        \Delta (\sd_{S_2}-\sd_{S_1 \times \{0\}}) \leq 0
        \quad 
        \text{on } S_1 \times (0,+\infty) \text{ in distributional sense.}
        \end{equation}
        \par 
        Consider the weight function $V \in \Lip(S_1 \times (0,+\infty))$ given by  $V=n^\delta \sd_{S_2}$. 
        We denote by $\Delta_w \bar{\sd}$ the weighted Laplacian of $\bar{\sd}$ w.r.t. $V$, i.e. the measure $\Delta_w \bar{\sd}:=\Delta \bar{\sd}-\nabla V \cdot \nabla \bar{\sd}$. Thanks to \eqref{E1} and the definition of $V$, it holds
        \begin{equation} \label{E3}
            \Delta_w \bar{\sd} \leq 0
            \quad 
        \text{on } S_1 \times (0,+\infty) \text{ in distributional sense.}
        \end{equation}
        Set $v:=(1+\bar{\sd})^{-1}$ and note that, by the chain rule, it holds in distributional sense on $S_1 \times (0,+\infty)$ 
        \begin{equation} \label{Efin2}
            \Delta_w v=-(1+\bar{\sd})^{-2}\Delta_w \bar{\sd}+2(1+\bar{\sd})^{-3}|\nabla  \bar{\sd}|^2 \geq 0. 
        \end{equation}
        
        Let $\phi \in \Lip_c(S_1 \times [0,+\infty))$ be a positive function.
        Consider the functions $\phi_\epsilon \in \Lip(S_1 \times [0,+\infty))$ defined by $\phi_\epsilon(x,t):=(1-\epsilon^{-1}t) \vee 0$.
        Integrating by parts, it holds
        \begin{align} \label{EFin1}
       & \int_{S_1 \times (0,+\infty)} (\nabla v \cdot \nabla \phi) \,  e^{-V} d \Vol \\ \nonumber
        &= \int_{S_1 \times (0,+\infty)} [\nabla v \cdot \nabla (\phi \phi_\epsilon)] \, e^{-V} d \Vol 
        + \int_{S_1 \times (0,+\infty)} [\nabla v \cdot \nabla (\phi (1-\phi_\epsilon)) ]\, e^{-V} d \Vol 
        \\ \nonumber
        & = \int_{S_1 \times (0,+\infty)} [\nabla v \cdot (\phi \nabla  \phi_\epsilon+\phi_\epsilon \nabla \phi)] \, e^{-V} d \Vol 
        -\int_{S_1 \times (0,+\infty)}  \phi (1-\phi_\epsilon) \, e^{-V} d \Delta_w v. 
        \end{align}
        Observe now that $\nabla \phi_\epsilon=-\epsilon^{-1}\nabla \sd_{S_1 \times \{0\}}1_{[0,\epsilon]}$, so that
        \begin{equation}
            \nabla v \cdot \nabla \phi_\epsilon=\epsilon^{-1}1_{[0,\epsilon]}(1+\bar{\sd})^{-2}\nabla \bar{\sd} \cdot \nabla \sd_{S_1 \times \{0\}} \leq 0.
        \end{equation}
        Hence,
         passing to the limit in \eqref{EFin1} as $\epsilon \to 0$ and taking into account \eqref{Efin2}, it holds
        \begin{equation} \label{E4}
            \int_{S_1 \times (0,+\infty)} (\nabla v \cdot \nabla \phi) \,  e^{-V} d \Vol \leq 0.
        \end{equation}
        Let $\psi \in \Lip_c(S_1 \times [0,+\infty))$ be a positive function. 
        
        We repeat a computation from \cite{Y} (see also \cite[Lemma $7.1$]{LiG}).
        By \eqref{E4}, it holds
        \begin{align*}
        0 &\geq \int_{S_1 \times (0,+\infty)} [\nabla (\psi^2v) \cdot \nabla v] \, e^{-V} d \Vol \\
        &=2 \int_{S_1 \times (0,+\infty)} \psi v (\nabla \psi \cdot \nabla v) \, e^{-V} d \Vol+\int_{S_1 \times (0,+\infty)} \psi^2 |\nabla v|^2 \, e^{-V} d \Vol,
        \end{align*}
        so that
        \begin{equation} \label{EE}
        \int_{S_1 \times (0,+\infty)} \psi^2 |\nabla v|^2 \, e^{-V} d \Vol \leq 2 \int_{S_1 \times (0,+\infty)}|\psi v 
        \nabla \psi \cdot \nabla v| \, e^{-V} d \Vol.
        \end{equation}
        By Young's inequality, we then obtain
        \begin{align*}
        2 \int_{S_1 \times (0,+\infty)}|\psi v 
        \nabla \psi \cdot \nabla v| \,  e^{-V} d \Vol
        \leq & \; \frac{1}{2} \int_{S_1 \times (0,+\infty)} \psi^2 |\nabla v|^2 \,  e^{-V} d \Vol \\
        &+8 \int_{S_1 \times (0,+\infty)} |\nabla \psi|^2 v^2 \,  e^{-V} d \Vol,
        \end{align*}
        which, combining with \eqref{EE}, gives
        \[
        \int_{S_1 \times (0,+\infty)} \psi^2 |\nabla v|^2 \, e^{-V} d \Vol \leq
        16 \int_{S_1 \times (0,+\infty)} |\nabla \psi|^2 v^2 \, e^{-V} d \Vol.
        \]
        Recalling the definition of $v$, this implies
        \begin{equation} \label{E5}
        \int_{S_1 \times (0,+\infty)} \psi^2 \frac{|\nabla \bar{\sd}|^2}{(1+\bar{\sd})^4} \, e^{-V}d\Vol \leq
        16 \int_{S_1 \times (0,+\infty)} |\nabla \psi|^2 \, e^{-V}d\Vol.
        \end{equation}
    Let $x \in S_1$ and $r>0$. Let $\phi_1 \in \Lip_c(S_1)$ be such that $\phi_1 \equiv 1$ in $B_r(x) \subset S_1$ and $\phi_2 \in \Lip_c([0,+\infty))$ be such that $\phi_2 \equiv 1$ in $[0,r)$, both taking values in $[0,1]$. Let $\psi \in \Lip_c(S_1 \times [0,+\infty))$ be defined by $\psi(y,t):=\phi_1(y)\phi_2(t)$. We denote by $|\cdot |_t$ and $\nabla^t \phi_1$ the norm and the gradient of $\phi_1$ w.r.t. the metric $e^{2\delta t}g_1$ on $S_1$, and we note that $|\nabla^t \phi_1|_t=e^{-\delta t} |\nabla^0 \phi_1|_0$. In particular, it holds
    \[
    |\nabla \psi|^2(y,t) \leq |\nabla^0 \phi_1|_0^2(y)+\phi_1^2(y)(\phi_2'(t))^2.
    \]
    Hence, with this choice of $\psi$, the r.h.s. in \eqref{E5} is bounded from above by
    \begin{equation} \label{E6}
    \begin{split}
       & 16
        \int_0^{+\infty}e^{-n^\delta t+\delta(n-1)t} \, \de t \int_{S_1} |\nabla^0 \phi_1|_0^2 \, d\Vol_{S_1}
       \\& \quad  + 16
        \int_{S_1} \phi_1^2 \, d\Vol_{S_1} \int_0^{+\infty} |\phi_2'|^2  e^{-n^\delta t+\delta(n-1)t} \, \de t.
        \end{split}
    \end{equation}
    Taking an infimum over all the admissible $\phi_2 \in \Lip_c([0,+\infty))$ such that $\phi_2 \equiv 1$ on $[0,r)$, we deduce
    \begin{equation} \label{E7}
    \int_{B_r(x) \times [0,r]} \frac{|\nabla \bar{\sd}|^2}{(1+\bar{\sd})^4} \, e^{-V}d\Vol \leq 16
        \int_0^{+\infty}e^{- t} \, \de t \int_{S_1} |\nabla^0 \phi_1|_0^2 \, d\Vol_{S_1}.
    \end{equation}
    Since $S_1$ is parabolic, taking the infimum of the r.h.s. over all the functions $\phi_1 \in \Lip_c(S_1)$ such that $\phi_1 \equiv 1$ on $B_r(x) \subset S_1$, we deduce $|\nabla \bar{\sd}| \equiv 0$ on $B_r(x) \times [0,r]$. The constancy of $\bar{\sd}$ follows from the arbitrariness of $r>0$, concluding the proof.
    \end{proof}
    \end{thm}
    
In the previous proof, the weight $V$ was needed to compensate for the excessive volume growth of the space $S_1 \times (0,+\infty)$ equipped with $$\de s^2=e^{2\delta t}g_1+\de t^2.$$ In the proof of Theorem \ref{CT2}, we consider $S_2 \times (0,+\infty)$ equipped with $$\de s^2=e^{-2 t}g_2+\de t^2,$$ $g_2$ being the metric on $S_2$. Hence, the volume growth will be easier to control even without a weight. 

Nevertheless, given $f \in C^\infty(S_2)$, its gradient w.r.t. $e^{-2 t}g_2$ increases as $t$ increases. Hence, also in the proof of Theorem \ref{CT2}, a weight is needed, this time to compensate for the previously described gradient growth. 
The next theorem implies Theorem \ref{CT2} from the Introduction. It is a slightly more general result that allows for multiple boundary components to satisfy $\mathrm{I\!I} \geq -1$.

\begin{thm}
     Let $(M^n,g)$ be a Riemannian manifold, with $\Ric_{M} \geq -(n-1)$. Let $\partial M$ be the disjoint union $\partial M=S_1\cup S_2\cup S_3$, where each $S_i$ is a union of connected components of $\partial M$. Assume that the following hold.
     \begin{enumerate}
        \item 
    $S_1$ is non-empty with $H_{S_1} \geq 1$.
         \item $S_2$ is non-empty, connected, parabolic, with $\mathrm{I\!I}_{S_2} \geq -1$, and $\Ric_{S_2} \geq 0$. \item 
     $S_3$ has $\mathrm{I\!I}_{S_3} \geq -1$, and $\Ric_{S_3} \geq 0$. 
     \end{enumerate}
    Then, $S_3=\emptyset$, and there exists $l>0$ such that $M$ is isometric to $S_1 \times [0,l]$ with
    the warped product metric $\de s^2 = e^{-2 t}g_1 + \de t^2$, where $g_1$ is the metric on $S_1$.
    \begin{proof}
         We argue as in Theorem \ref{CT1body}.
        Let $$S:=S_2 \cup S_3.$$ Consider the manifold $S \times (0,+\infty)$ equipped with the warped product metric $$\tilde{g} = e^{-2t}g_S + \de t^2,$$ where $g_S$ is the metric induced by $M$ on $S$. By standard computations, the second fundamental form of $S \times \{0\}$ in $S \times (0,+\infty)$ satisfies
        $
        \mathrm{I\!I}_{S \times \{0\}}=1
        $. Moreover, since $\Ric_{S} \geq 0$, Gauss' equations imply that $$\Ric_{S \times (0,+\infty)} \geq -(n-1).$$
        
        In a similar fashion, we denote by $g_2$ the restriction of $g_S$ to $S_2$ and we consider $S_2 \times (0,+\infty)$, which is a connected component of $S\times(0,+\infty)$. We denote by $\ssf{Vol}$ the Riemannian volume in $S_2 \times (0,+\infty)$ w.r.t. the metric $\tilde{g}$ and by $\ssf{Vol}_{S_2}$ the Riemannian volume on $S_2$ w.r.t. the metric $g_2$.
        Let $\sd_{S_2 \times\{0\}}$ be the distance from $S_2 \times \{0\}$ in $S_2 \times (0,+\infty)$.
        It holds
        $$
        \Delta \sd_{S_2 \times \{0\}}=-(n-1) \quad \text{on } S_2 \times (0,+\infty).$$
        Consider the metric space $(\X,\sd)$ obtained by gluing $M$ and $S \times (0,+\infty)$ along the two boundaries isometric to $S$. On this metric space, we consider the distance function from $S_1$ and we denote it by $\sd_{S_1}$.
        In particular, this allows us to define $\sd_{S_1}$ on $S_2 \times (0,+\infty)$ (seen as a subset of the previously defined metric space). 
        We first show that the function $$\bar{\sd}:=\sd_{S_1}-\sd_{S_2 \times \{0\}}$$ is constant on $S_2 \times (0,+\infty)$. 
        
        By Lemma \ref{l|Laplacian in the glued space}, it holds $\Delta \sd_{S_1} \leq -(n-1)$ in distributional sense on $S_2 \times (0,+\infty)$. Hence, 
        \begin{equation}  \label{E11}
        \Delta (\sd_{S_1}-\sd_{S_2 \times \{0\}}) \leq 0
        \quad 
        \text{on } S_2 \times (0,+\infty) \text{ in distributional sense.}
        \end{equation} 
       Consider the weight function $V \in \Lip(S_2 \times (0,+\infty))$ given by  $$V=2 \sd_{S_1}.$$ 
        Let $\phi \in \Lip_c(S_2 \times [0,+\infty))$ be a positive function. 
        Arguing as in Theorem \ref{CT1body} (derivation of \eqref{E4}), one obtains, setting $v:=(1+\bar{\sd})^{-1}$,
        \begin{equation} \label{E8}
        \int_{S_2 \times (0,+\infty)} (\nabla v \cdot \nabla \phi) \, e^{-V} d \Vol \leq 0. 
        \end{equation}
        If $\psi \in \Lip_c(S_2 \times [0,+\infty))$ is a positive function, arguing again as in Theorem \ref{CT1body} (derivation of \eqref{E5}), one has 
        \begin{equation} \label{E9}
        \int_{S_2 \times (0,+\infty)} \psi^2 \frac{|\nabla \bar{\sd}|^2}{(1+\bar{\sd})^4} \, e^{-V} d\Vol \leq
        16 \int_{S_2 \times (0,+\infty)} |\nabla \psi|^2 \, e^{-V} d\Vol.
        \end{equation}
    Let $x \in S_2$ and $r>0$. Let $\phi_1 \in \Lip_c(S_2)$ be such that $\phi_1 \equiv 1$ in $B_r(x)$ and $\phi_2 \in \Lip_c([0,+\infty))$ be such that $\phi_2 \equiv 1$ in $[0,r)$, both taking values in $[0,1]$. Let $\psi \in \Lip_c(S_2 \times [0,+\infty))$ be defined by $\psi(y,t):=\phi_1(y)\phi_2(t)$. 
    We denote by $|\cdot|_t$ and $\nabla^t \phi_1$ the norm and the gradient of $\phi_1$ w.r.t. the metric $e^{-2t}g_2$ on $S_2$, and we note that $|\nabla^t \phi_1|_t=e^{t} |\nabla^0 \phi_1|_0$. In particular, it holds
    \[
    |\nabla \psi|^2(y,t) \leq e^{2t}|\nabla^0 \phi_1|_0^2(y)+\phi_1^2(y)(\phi_2'(t))^2.
    \]
    With this choice of $\psi$, the r.h.s. in \eqref{E9} is bounded from above by
    \begin{equation} \label{E10}
    \begin{split}
        16&
        \int_0^{+\infty}e^{-(n-1)t} \, \de t \int_{S_2} |\nabla^0 \phi_1|_0^2 \, d\Vol_{S_2} \\
        &+ 16
        \int_{S_2} \phi_1^2 \, d\Vol_{S_2} \int_0^{+\infty} (\phi_2')^2  e^{-(n-1)t} \, \de t.
        \end{split}
    \end{equation}
    We now conclude that $\bar{\sd}$ is constant on $S_2 \times (0,+\infty)$ as in Theorem \ref{CT1body}.
    Hence, $\sd_{S_1}$ is constant on $S_2$ with constant value $D>0$. 
    
    Let now $\sd_{S_2}$ be the distance from $S_2$ in $\X \setminus S_2 \times [0,+\infty)$. By Lemma \ref{l|Laplacian in the glued space}, it holds $\Delta (\sd_{S_1}+\sd_{S_2} )\leq 0$ in the $D$-neighbourhood of $S_2$ in $\X \setminus S_2 \times [0,+\infty)$. 
    We claim that the set 
    \[
    C:=\{x \in \X \setminus (S_1 \cup S_2 \times [0,+\infty)): \sd_{S_1}+\sd_{S_2}=D\}
    \]
    is open and closed in $\X \setminus (S_1 \cup S_2 \times [0,+\infty))$. 
    Closedness is trivial. Let $p \in C$. We show that there exists a neighbourhood of $p$ in $\X$ contained in $C$. 
    The set $C$ is contained in the $D$-neighbourhood of $S_2$, while $\sd_{S_1}+\sd_{S_2} \geq D$ on $\X$. Hence, by the maximum principle, $\sd_{S_1}+\sd_{S_2} = D$ in a neighbourhood of $p$. Since $C$ is open and closed in $\X \setminus (S_1 \cup S_2 \times [0,+\infty))$, it coincides with $\X \setminus (S_1 \cup S_2 \times [0,+\infty))$ itself. Hence, $S_3=\emptyset$, since otherwise $\sd_{S_1}+\sd_{S_2}$ would attain arbitrarily large values on $S_3 \times (0,+\infty) \subset \X \setminus (S_1 \cup S_2 \times [0,+\infty))$.
    The conclusion then follows from Lemma \ref{CL2}.
    \end{proof}
    \end{thm}

\section{Applications}

Corollaries \ref{CC2} and \ref{CC3:intro} will be consequences of Theorem \ref{CT4} below, which relates the problem of deciding whether a manifold with disconnected mean-convex boundary is a product with the following well studied question:
\begin{itemize}
    \item Given a manifold $(M,g)$, under which conditions is a stable minimal hypersurface $\Sigma \subset M$ totally geodesic?
\end{itemize}

\begin{thm} \label{CT4}
    Let $n \leq 7$ and let $(M^n,g)$ be a manifold with non-negative Ricci curvature and mean-convex disconnected boundary. Assume that any properly embedded, two-sided, stable minimal hypersurface $\Sigma \subset M$ is parabolic and has $\Ric_{\Sigma} \geq 0$. Then, $(M,g)$ splits isometrically as $\Sigma \times [0,l]$, for some manifold $(\Sigma,g')$ and $l>0$.
    \begin{proof}
        If one boundary component $\Sigma$ of $\partial M$ is minimal and stable, then by assumption it is parabolic and it has $\Ric_{\Sigma} \geq 0$. By the second variation formula, $\Sigma$ is also totally geodesic. Hence, the claim follows by Theorem \ref{CT1}.

        If no connected component of $\partial M$ is minimal and stable, we can consider the set $E \subset M$ given by Lemma \ref{L|yau} relative to a connected component $\Sigma$ of $\partial M$. Since $n \leq 7$, $E$ is a smooth set. Since $E$ minimizes the perimeter w.r.t. inner competitors, its boundary $\partial E$ is a properly embedded, two-sided, stable minimal hypersurface. By assumption, each connected component of $\partial E$ is parabolic and it has $\Ric_{\partial E} \geq 0$. By Theorem \ref{CT1}, it follows that $E=\Sigma \times [0,l]$ isometrically for some $l>0$. Applying the same argument to $M \setminus E$, we conclude.
    \end{proof}
\end{thm}

\begin{remark}[On previous related results]
The idea of using stable minimal hypersurfaces to deduce existence of a product structure in a manifold with non-negative Ricci curvature is classical, and it has been used, for instance, in \cite{AndersonRodriguez,Anderson,GangLiu}. In the aforementioned results, a key step to prove existence of a product structure is to build a \emph{foliation} of stable minimal hypersurfaces in the manifold in question. In \cite{AndersonRodriguez,Anderson}, this is possible thanks to a uniform bound on the sectional curvature. In \cite{GangLiu}, the existence of such a foliation follows from an assumption on the fundamental group. Theorem \ref{CT4}, being a consequence of Theorem \ref{CT1}, does not require the construction of such a foliation. 

Finally, we mention that there
are many other splitting results that can be obtained with different lower curvature  bounds, and that involve the existence (or the construction) of either minimal hypersurfaces or $\mu$-bubbles (see \cite{zGromovbubbles1,GromovLecturesScalar}). Among many others, \cite{zbMATH07063922,zbMATH06664766} obtain splitting results for $3$-manifolds with non-negative scalar curvature,
while splitting results for manifolds with Ricci curvature bounded from below in spectral sense are obtained in \cite{arXiv:2412.12707, arXiv:2412.12631,zbMATH08159607}.
\end{remark}

The next result corresponds to Corollary \ref{CC2}.
\begin{thm}
    Let $(M^3,g)$ be a Riemannian manifold with non-negative Ricci curvature and disconnected mean-convex boundary. Then, $M=\Sigma \times [0,l]$ isometrically, for some manifold $(\Sigma,g')$ and $l>0$.
\begin{proof}
By \cite{SchoenYauStability}, any properly embedded, two-sided, stable minimal hypersurface $\Sigma \subset M$ is totally geodesic, and the Ricci curvature of $M$ in the normal direction to $\Sigma$ vanishes. In particular, $\Sigma$ has non-negative sectional curvature. Hence, by Bishop-Gromov's inequality, $\Sigma$ is parabolic.
The statement then follows by Theorem \ref{CT4}.
\end{proof}
\end{thm}

Before turning our attention to Corollary \ref{CC3:intro}, we briefly discuss the notions of $(n-2)$-Ricci curvature and weakly bounded geometry. 

\begin{definition}[Lower bounds on $(n-2)$-Ricci curvature] \label{D2}
    A manifold $(M^n,g)$ is said to have non-negative $(n-2)$-Ricci curvature (denoted $\Ric_{n-2} \geq 0$) if for every $p \in M$ and every collection of orthonormal vectors $$\{e_1, \cdots, e_{n-1}\} \subset T_pM,$$ it holds
    \[
    \sum_{i=1}^{n-2} \ssf{Sec}(e_{n-1},e_i) \geq 0.
    \]
\end{definition}

The condition $\Ric_{n-2} \geq 0$ is an intermediate condition between $\mathsf{Sec} \geq 0$ and $\Ric \geq 0$, in the sense that $\mathsf{Sec} \geq 0 \Rightarrow \Ric_{n-2} \geq 0 \Rightarrow \Ric \geq 0$. 
An important feature of manifolds $(M^n,g)$ with $\Ric_{n-2} \geq 0$ is the following: if $\Sigma^{n-1} \subset M^n$ is a totally geodesic hypersurface and $\Ric_{n-2} \geq 0$ on $M^n$, then $\Ric_\Sigma \geq 0$ on $\Sigma$. This is an immediate consequence of the Gauss equation since, having fixed an orthonormal basis $\{e_1,\dots,e_{n-1}\}$ of $T_p\Sigma$, we have
\begin{equation} \label{EFin10}
    \Ric_\Sigma(e_{n-1},e_{n-1}):=\sum_{i=1}^{n-2} \ssf{Sec}_\Sigma(e_{n-1},e_i)=\sum_{i=1}^{n-2} \ssf{Sec}_M(e_{n-1},e_i) \geq 0.
\end{equation}


We now recall the definition of weakly bounded geometry (after \cite[Definition 2.3]{chodosh2024completestableminimalhypersurfaces}).
 \begin{definition}[Weakly bounded geometry] \label{D1}
A manifold $(M^n,g)$ without boundary is said to have weakly bounded geometry if there exists $Q>0$ and $\alpha \in (0,1)$ satisfying the following. For every $x \in M$, there is an open set $U_x \ni x$ and a $C^{2,\alpha}$ \emph{local} diffeomorphism $\phi:B^{\bb{R}^n}_{Q^{-1}}(0) \to U_x$ such that
\begin{enumerate}
    \item $e^{-2Q}g_{eu} \leq \phi^* g \leq e^{2Q}g_{eu}$ as bilinear forms.
    \item $\|\partial_k \phi^*g_{i,j}\|_{C^{0,\alpha}} \leq Q$, for every $i,j,k \in \{1,\dots,n\}$.
\end{enumerate}
\end{definition}

We recall that a manifold $(M^n,g)$ is said to have $m$-bounded geometry (see \cite[Section 3]{cheegergromovtaylor})
if there exists $Q>0$ satisfying the following. For every $x \in M$, there is a geodesic ball $B_{Q^{-1}}(x)$ such that the exponential map $\mathrm{exp}_x:B^{\bb{R}^n}_{Q^{-1}}(0) \to B_{Q^{-1}}(x)$ is a diffeomorphism and 
\begin{enumerate}
    \item $e^{-2Q}g_{eu} \leq \mathrm{exp}_x^* g \leq e^{2Q}g_{eu}$ as bilinear forms.
    \item $\| \mathrm{exp}_x^* g_{i,j}\|_{C^m} \leq Q$, for every $i,j,k \in \{1,\dots,n\}$.
\end{enumerate}

The main difference between weakly bounded geometry according to Definition \ref{D1} and $2$-bounded geometry, besides the $C^{1,\alpha}$ control on the pull back of the metric being weaker than $C^2$ control, is that Definition \ref{D1} does not force a uniform positive lower bound on the volume of balls or on the injectivity radius. To see this, it is enough to consider a hyperbolic cusp, as observed in \cite{chodosh2024completestableminimalhypersurfaces}. For more background on the notion of weakly bounded geometry, we refer to the original reference \cite{chodosh2024completestableminimalhypersurfaces}.

\begin{remark}[Geometric motivations for $\Ric_{n-2} \geq 0$, $\mathrm{Scal} \geq 1$, and weakly bounded geometry] \label{RGromovconj}
    Recently, there has been growing interest in the study of manifolds $(M^n,g)$ with $\Ric \geq 0$ and $\mathrm{Scal} \geq 1$ (see, for instance, \cite{ScalmeetsRicci,wanggrowth,Chodoshvolumethreemanifolds}). The main conjecture on these manifolds is due to Gromov \cite{LargeGromov}, and states that one should have
    \begin{equation}
        \limsup_{r \to + \infty} \frac{\mathsf{Vol}(B_r(p))}{r^{n-2}} < + \infty.
    \end{equation}
    If the conjecture holds, one consequence would be the following: $4$-manifolds with $\Ric \geq 0$, $\mathrm{Scal} \geq 1$, and $b_1 \neq 0$ have universal covering that splits off a line isometrically (compare with Theorem \ref{T|Anderson}). This would be true because a $4$-manifold as above has at most linear volume growth at infinity by \cite{Anderson}, so that its universal covering splits off a line by \cite{HuangHuang}.

    Manifolds with $\Ric_{n-2} \geq 0$, $\mathrm{Scal} \geq 1$, and weakly bounded geometry provide a convenient subclass in which to study the geometric consequences of combined lower bounds on Ricci and Scalar curvature.
    
    In dimension $4$, these stronger assumptions imply that stable minimal hypersurfaces are totally geodesic by \cite{chodosh2024completestableminimalhypersurfaces} (compare with \cite[Example 1.9]{chodosh2024completestableminimalhypersurfaces}). This fact is key in proving  is Corollary \ref{CC3:intro} and Theorem \ref{T|Anderson}. We remark that the assumption of weakly bounded geometry is only needed to apply the results from \cite{chodosh2024completestableminimalhypersurfaces}, and is otherwise not used in the paper.
\end{remark}

The next result corresponds to Corollary \ref{CC3:intro}.

\begin{thm}
    Let $(M^4,g)$ be a Riemannian manifold with $\Ric_2 \geq 0$, scalar curvature $\geq 1$, and weakly bounded geometry. Let $N^4 \subset M^4$ be a smooth submanifold with mean-convex disconnected boundary. Then, $N=\Sigma \times [0,l]$ isometrically, for some manifold $(\Sigma,g')$ and $l>0$.
    \begin{proof}
        By \cite{chodosh2024completestableminimalhypersurfaces}, any properly embedded, two-sided, stable minimal hypersurface $\Sigma \subset N$ is parabolic and totally geodesic.
The statement then follows by Theorem \ref{CT4}.
    \end{proof}
\end{thm}

The next result corresponds to Theorem \ref{T|Anderson}.

\begin{thm}
    Let $(M^n,g)$ be a manifold without boundary such that
    either
    \begin{itemize}
    \item  $n=3$, and $\Ric \geq 0$,
    \end{itemize}
    or
    \begin{itemize}
    \item $n=4$, $\Ric_2 \geq 0$, $\mathrm{Scal} \geq 1$, and $M$ has weakly bounded geometry.
    \end{itemize}
    If $b_1(M) \neq 0$, then the universal covering of $M$ splits a line isometrically.
    \begin{proof}
    Modulo passing to a covering of the manifold, we may assume that $\pi_1(M) \cong \bb{Z}$.
        Let $\gamma \subset M$ be a smooth loop generating $\pi_1(M)$.
           By the argument of  \cite[p.217, proof of Lemma 3]{SchoenYauStability}, there exists a properly embedded stable minimal hypersurface $\Sigma \subset M$ intersecting (a perturbation of) $\gamma$ transversally with algebraic intersection number $1$. By \cite{SchoenYauStability} and \cite{chodosh2024completestableminimalhypersurfaces}, $\Sigma$ is totally geodesic and satisfies $\Ric(\nu,\nu) \equiv 0$, $\nu$ being the normal to $\Sigma$. 
           \smallskip
           
           \textbf{Step 1}: $\Sigma$ is
           parabolic, and satisfies $\Ric_\Sigma \geq 0$. 
           
           Let us consider first the case where $n=3$ and $\Ric_M \geq 0$. In this case, the Gauss equations together with the fact that $\Ric(\nu,\nu) \equiv 0$ imply that $\Sigma$ has non-negative scalar curvature. Since $\Sigma$ has dimension $2$, then it has non-negative sectional curvature and it is parabolic by Bishop-Gromov's inequality.

           Let us consider now the case where $n=4$, and $M$ has $\Ric_{2} \geq 0$, $\mathrm{Scal} \geq 1$ and weakly bounded geometry. By the Gauss equations, the assumption  $\Ric_{2} \geq 0$ implies that  $\Ric_\Sigma \geq 0$ (see the discussion around \eqref{EFin10}).
           Moreover, combining  the Gauss equations, the fact that $\Ric(\nu,\nu) \equiv 0$, and $\mathrm{Scal} \geq 1$, it follows that $\Sigma$ has itself scalar curvature greater than $1$. In particular, by \cite{wanggrowth}, $\Sigma$ has at most linear volume growth so that it is parabolic as claimed.
\smallskip

           \textbf{Step 2:}
           The image of the map $i_\#:\pi_1(\Sigma) \to \pi_1(M)$ induced by the inclusion is trivial. 
           
           Indeed, if this were not the case, then $i_{\#}(\pi_1(\Sigma)) \subset \pi_1(M)$, being a subgroup of $\pi_1(M) \cong \bb{Z}$, would contain a copy of $\mathrm{span}(\gamma^k)$ for some $k \in \bb{N}$. Hence, $\gamma^k$ can be homotopically deformed into a loop fully contained in $\Sigma$, which can be then deformed into a loop which is disjoint from $\Sigma$ (moving away from $\Sigma$ in its tubular neighbourhood). Hence, the intersection number of $\gamma^k$ with $\Sigma$ is zero. At the same time, since the intersection number of $\gamma$ and $\Sigma$ is one, it follows that the intersection number of $\Sigma$ with $\gamma^k$ is $k$, a contradiction.
\smallskip

     \textbf{Step 3:}  conclusion of the proof. 
     
     Let $p:\tilde{M} \to M$ be the universal covering of $M$, and let $q \in \Sigma \cap \gamma$. Let $\Sigma_1, \Sigma_2$ be lifts of $\Sigma$ passing through distinct points $q_1,q_2 \in \tilde{M} \cap p^{-1}(q)$. Observe that these exist by Step $2$ and the standard lifting criterion. These two hypersurfaces are disjoint, since otherwise we would lose uniqueness of the lift at the points that they share. Consider an $n$-manifold with boundary $V \subset \tilde{M}$ such that $\partial V=\Sigma_1 \cup \Sigma_2$ topologically (such $V$ exists since $\tilde{M}$ is simply connected, and $\Sigma_1,\Sigma_2$ are properly embedded codimension-1 submanifolds without boundary). 
     By Corollaries \ref{CC2} and \ref{CC3:intro}, it holds $V \cong \Sigma \times [0,l]$ isometrically for some $l>0$.
     Iterating this procedure for every pair of points $q_1,q_2 \in \tilde{M} \cap p^{-1}(q)$, we conclude. 
    \end{proof}
\end{thm}



\begin{remark} [On higher dimensional versions of Corollaries \ref{CC2}, \ref{CC3:intro}, and Theorem \ref{T|Anderson}] Both Corollary \ref{CC2} and Theorem \ref{T|Anderson} (when considered under the assumptions of Corollary \ref{CC2}) fail in higher dimensions. Indeed, there is a metric with $\Ric>0$ on $S^1 \times \bb{R}^3$ constructed by Nabonnand in \cite{zbMATH03709169}. By the procedure used in the proof of Theorem \ref{T|Anderson}, in the universal covering of $S^1 \times \bb{R}^3$, there is a region $V$ whose boundary is disconnected and consists of minimal hypersurfaces. Nevertheless, this region has $\Ric>0$. In particular,  it cannot split.

On the other hand, Corollary \ref{CC3:intro} and Theorem \ref{T|Anderson} (when considered under the assumptions of Corollary \ref{CC3:intro}) are expected to hold in dimension $5$. This would follow (by adapting the arguments above) if one could show that, in a $5$-manifold with $\Ric_3 \geq 0$, $\mathrm{Scal} \geq 1$, and weakly bounded geometry, any boundary of a perimeter-minimizing set is totally geodesic.

There are several reasons to believe that this is the case (see \cite[Remark 1.8]{chodosh2024completestableminimalhypersurfaces}). One indication is that, assuming Gromov’s volume growth conjecture mentioned in Remark \ref{RGromovconj}, such manifolds have at most cubic volume growth at infinity. It would then follow that the boundaries of perimeter-minimizing sets (endowed with their intrinsic metric) are parabolic. The desired conclusion would then follow from the second variation formula combined with the logarithmic cut-off argument.

We note that this strategy would only require $\Ric_{3} \geq 0$ and $\mathrm{Scal} \geq 1$, without the need for the weakly bounded geometry assumption.
\end{remark}

Finally, the next result corresponds to Corollary \ref{CC1}.

\begin{thm}
    Let $\delta \in \{0,1\}$. Let $(M^n,g)$ be a parabolic manifold with $\Ric_M \geq 0$ and let $M \times \bb{R}$ be equipped with $\de s^2 = e^{-2 t \delta}g + \de t^2$. If $E \subset M \times (0,+\infty)$ is a smooth closed set with connected boundary and mean curvature $H_{\partial E} \leq \delta$, then $E=M \times [a,+\infty)$, for some $a > 0$.
    If $E \subset M \times (-\infty ,0)$ is a smooth closed set with connected boundary and mean curvature $H_{\partial E} \leq -\delta$, then $E=M \times (-\infty,a]$, for some $a < 0$.
\begin{proof}
Let $E \subset M \times (0,+\infty)$ be with mean curvature satisfying $H_{\partial E} \leq \delta$. Consider the manifold 
$$N:=M \times (0,+\infty) \setminus E.$$ Then, $\partial N=\partial E \cup M \times \{0\}$. Since the second fundamental form of $M \times \{0\}$ in $M \times (0,+\infty)$ is equal to $\delta$, the result follows from Theorem \ref{CT1} 9f $\delta=0$ and Theorem \ref{CT1*} if $\delta=1$.
The case when $E \subset M \times (-\infty,0)$ has mean curvature satisfying $H_{\partial E} \leq -\delta$ follows analogously, applying Theorem \ref{CT2} instead of Theorem \ref{CT1*}.
\end{proof}
\end{thm}

\section{Conclusive remarks} \label{S5}

In this final section, we give a partial answer to the question posed in Remark \ref{R2}. In Theorems \ref{T|S5T1} and \ref{T|S5T2} below, we prove variants of Theorem \ref{CT1} where the assumptions on the second fundamental form and the Ricci curvature of the parabolic boundary component are removed. Nevertheless, we compensate by strengthening the curvature assumptions on the ambient space ($\Ric_{n-2} \geq 0$ instead of $\Ric \geq 0$), and by requiring different additional conditions on a boundary component.
    
    \begin{thm} \label{T|S5T1}
    Let $(M^n,g)$ be a  manifold with $\Ric_{n-2} \geq 0$ and disconnected mean-convex boundary. Let $\Sigma \subset \partial M$ be a boundary component which is minimal, stable, and parabolic.
    Then, $M=\Sigma \times [0,l]$ isometrically for some $l>0$.
    \begin{proof}
        Since $\Sigma$ is minimal, stable, and parabolic, then it is totally geodesic by the second variation formula of the area. Since $M$ has $\Ric_{n-2} \geq 0$, we then obtain that $\Ric_\Sigma \geq 0$. Hence, the statement follows by Theorem \ref{CT1}.
    \end{proof}
\end{thm}

 \begin{thm} \label{T|S5T2}
    Let $(M^n,g,p)$ be a pointed manifold with $\Ric_{n-2} \geq 0$ and disconnected mean-convex boundary. Let $\Sigma \subset \partial M$ be a boundary component with
    \begin{equation} \label{E|volgrowth}
        \int_{1}^\infty \frac{t}{\aH^{n-1}(\partial B_t(p))+\aH^{n-1}(B_t(p) \cap \Sigma)} \, \de t=\infty.
    \end{equation}
    Then, $M=\Sigma \times [0,l]$ isometrically for some $l>0$.
    \begin{proof}
   We apply Lemma \ref{L|yau}, and we observe that if item \ref{item1Yau} of Lemma \ref{L|yau} is satisfied, taking into account our volume growth assumption, we find a connected component of $\partial M$ which is minimal, stable, and parabolic (due to the assumption \eqref{E|volgrowth}). In this case the conclusion follows by Theorem \ref{T|S5T1}.
        
        Hence, consider the set $E \subset M$ provided by item \ref{item2Yau} of Lemma \ref{L|yau} relative to $\Sigma$. By \eqref{E|volgrowth}, it holds
        \begin{equation} \label{E|par}
        \int_1^{\infty} \frac{t}{\aH^{n-1}(\partial E \cap B_t(p))} \, dt=+\infty.
        \end{equation}
        By repeating the proof of \cite[Theorem 3.1]{cucinotta2025perimeterminimizingsetsmanifolds} (cf.\ \cite[Theorem 2.1]{Anderson}), and using that $E$ minimizes the perimeter in $M \setminus \partial M$, it follows that $\partial E$ is smooth and totally geodesic. Since $\Ric_{n-2} \geq 0$, we deduce $\Ric_{\partial E} \geq 0$. Moreover, \eqref{E|par} implies that each connected component of $\partial E$ is parabolic. Hence, $E$ is a manifold with boundary, whose boundary is the union of its topological boundary $\partial E \subset M$ in $M$, and $\Sigma$. By Theorem \ref{CT1}, $E=\Sigma \times [0,l]$ isometrically for some $l>0$. Applying the same argument to $M \setminus E$, we obtain the statement.
    \end{proof}
\end{thm}

    


\renewcommand*{\bibfont}{\normalfont\small}

\small{
\printbibliography
}

\end{document}